\RequirePackage{fix-cm}
\documentclass[smallextended]{svjour3}       
\smartqed  
\usepackage{graphicx}

 \usepackage{amssymb}      
\usepackage{color}
\newcommand{\bu}{\mathbf u}
\newcommand{\bv}{\mathbf v}

\newcommand{\bU}{\mathbf U}
\newcommand{\bV}{\mathbf V}

\title{Using BDF schemes in the temporal integration of POD-ROM methods}
\author{Bosco
Garc\'{\i}a-Archilla\thanks{Departamento de Matem\'atica Aplicada
II, Universidad de Sevilla, Sevilla, Spain. Research is supported by
Spanish MCINYU under grants PID2021-123200NB-I00 and PID2022-136550NB-I00. (bosco@esi.us.es)}
\and Alicia Garc\'{\i}a-Mascaraque\thanks{Departamento de
Matem\'aticas, Universidad Aut\'onoma de Madrid, Spain. Research is supported
by predoctoral contract associated to PID2022-136550NB-I00 supported by MCIN/AEI/10.13039/501100011033 and FSE+ (alicia.garcia-mascaraque@uam.es)}
  \and Julia Novo\thanks{Departamento de
Matem\'aticas, Universidad Aut\'onoma de Madrid, Spain. Research is supported
by Spanish MINECO
under grant PID2022-136550NB-I00. (julia.novo@uam.es)}}
\date{\today}

\begin{document}
\maketitle
\abstract{In this paper we {consider} the numerical approximation of a semilinear reaction-diffusion model problem
by means of reduced-order methods (ROMs) based on proper orthogonal decomposition (POD). We focus on the time integration of the fully discrete reduced-order model. Most of the analysis in the literature has been carried out for the implicit Euler method
as time integrator. We integrate in time the reduced-order model with the BDF-q time stepping ($1\le q\le 5$) and prove optimal rate of convergence of order $q$ in time. Our set of snapshots is obtained from finite element approximations to the original model problem computed at different times. These finite element approximations can be obtained with any time integrator. 
The POD method is based on first order difference quotients of
the snapshots. The reason for doing this is twofold. On the one hand, the use of difference quotients allow us to provide pointwise-in-time error bounds. On the other, the use of difference quotients is essential to get the expected rate $q$ in time since
we  apply that the BDF-q time stepping, $1\le q\le 5$, can be written as a linear combination of first order difference quotients.}
\bigskip

\section{Introduction}
It is well known that using reduced-order models one can drastically reduce the computational cost of numerical simulations while
keeping enough accuracy in the desired approximations. In this paper we consider the so called proper orthogonal decomposition (POD) methods, \cite{ku-Vol}. 

The literature concerning POD methods has increased considerably in {\color{red}recent} years although \textcolor{red}{little} attention has been paid to the error coming from the temporal discretization. Most of the papers concerning error analysis apply the implicit Euler method as time integrator for the reduced-order model or at most second-order methods in time. However, in practice, higher order time integrators, including higher order BDF schemes are used. The use of higher order methods in time allow to take smaller time steps and, as a consequence, can also reduce the computational effort to get the final approximation. \textcolor{red}{The error analysis of the present paper has two main difficulties. The first one is the analysis of high order in time BDF-q methods for nonlinear problems. This difficulty is not associated to a reduced-order method, since it also appears when analyzing the standard Galerkin method. To solve this problem, we follow the error analysis in the recent paper \cite{bdf_stokes}. The $G$-stability results of BDF schemes for $q\ge 3$ and the clear application of those results to the error analysis of Galerkin methods shown in \cite{bdf_stokes} has been a key point to extend those techniques to the nonlinear problem we handle. The second difficulty is intrinsic to the reduced POD method. In order to bound the temporal error of the method, we need to estimate the projection error of the approximation to the time derivative in the BDF-q method (here, projection refers to the orthogonal projection onto the educed-order space) when applied to the snapshots. In a standard Galerkin method one has to bound only the error in the BDF-q approximation to the time derivative, which is known to have order $q$. The key point to bound the extra term is the inclusion of first order finite differences in the set of snapshots.}

 Since the literature about POD methods is quite 
extensive we will mention in this paper only those references that are related with the analysis we carry out in the present paper. \textcolor{red}{We refer the reader to the survey \cite{il_sur} and the references therein.} 

In the note \cite{letter}, the authors prove that second-order error bounds in time can be obtained for POD methods using the BDF2 scheme in time to approach the heat equation. The authors include difference quotients in the set of snapshots. The aim of the present paper, is to extend some of the ideas of \cite{letter} to get order $q$ in case of using BDF-q schemes for $q\le 5$ for semilinear reaction-diffusion problems. Then, we extend the ideas in \cite{letter} to higher order methods and to more complicated model problems. As in \cite{letter}, the use of difference quotients in the POD method is necessary to get pointwise-in-time error estimates. Moreover, the use of difference quotients allow us to get optimal bounds in time since we prove that any BDF-q method can be written as a linear combination of first order difference quotients (this idea was applied in \cite{letter} for the BDF2 scheme).

For the error analysis of this paper, we consider $H_0^1$-projections for the POD method since, as shown in \cite{rubino_etal}, better bounds in terms of the tail of the eigenvalues can be proved than in the case of $L^2$-projections. By tail of the eigenvalues we mean the quantity
\begin{equation}
\label{the_sigma_r}
\Sigma_r^2\equiv \sum_{k>r} \lambda_k
\end{equation}
where $\lambda_1\ge \lambda_2 \ge \ldots \ge \lambda_{d_r}>0$ are nonzero eigenvalues of the correlation matrix of the data set in the POD method. The reason is that
when projecting onto $L^2$ they are $H_0^1$ errors (concerning projections onto the reduced-order space) that have to be estimated 
and that are optimal in $L^2$, by definition, but not in $H_0^1$. In any case, the error analysis of the present paper also works for
$L^2$ projections. In that case, the analysis is simpler, not optimal in the sense of \cite{rubino_etal}, as any other
error analysis for $L^2$ projections, and the idea of writing any BDF-q scheme as linear combination of first order finite differences is not needed.

Apart from the idea of writing BDF-q methods as linear combination of first order schemes the key point for the error analysis of this paper is the analysis appearing in the recent reference \cite{bdf_stokes}. In \cite{bdf_stokes} the authors analyze time-stepping methods for BDF schemes of order 1 to 6 for the approximation of the transient Stokes equations with both inf-sup stable and symmetric pressure stabilized mixed finite element methods. As explained in \cite{bdf_stokes}, there are some results concerning
the $G$-stability of BDF schemes that allow the use of quite standard estimates that have made it easier to carry out the error analysis of
the POD method in this paper. The analysis of the BDF-6 method, presented in \cite{bdf_stokes}, could also be adapted to get order $6$ in time in the present paper. Nevertheless, this error analysis is harder and has to be done separately. 

For the error analysis of the POD method we also follow the ideas in \cite{temporal_nos} where a continuous-in-time POD method is analyzed. As in \cite{temporal_nos}, we assume our set of snapshots is based on  continuous-in-time finite element approximations to the model problem at different times. In practice, any time integrator can be applied to compute the snapshots, not necessarily the same used for the time integration of the POD method. Since the computation of snapshots is done in the offline stage, we can assume the snapshots are computed with enough accuracy so that the temporal errors coming from the time discretization are negligible. 
Although the error analysis of the present paper could be seen, in some sense, as the fully discrete case of the error
analysis in \cite{temporal_nos}, the ideas of \cite{temporal_nos} go further. In particular, interpolation techniques are applied in \cite{temporal_nos} to prove error estimations that are valid for values of the time variable different that those in the set of snapshots. 
As in \cite{temporal_nos}, the POD method is based on difference quotients of the snapshots both to get pointwise estimates and to be able to bound the errors coming from the temporal discretization.

Finally, we want to mention that one could also use the standard POD method based on the snapshots instead of its difference quotients and follow the error analysis in \cite{Bosco_pointwise}. In that paper, it is shown that a rate of convergence 
for a standard POD method not including difference quotients is as close as that of methods including difference quotients as the smoothness of the function providing the snapshots allow. In \cite{Bosco_pointwise} error bounds for the projections onto the POD space of the difference quotients are also provided, which would be needed to get the temporal estimations corresponding to the BDF-q scheme.

The outline of the paper is as follows, in Section 2 we introduce some preliminaries and notations. In Section 3 we state
some of the properties of the BDF methods appearing in \cite{bdf_stokes} that are needed to bound the time discretization error. In Section 4 we introduce some ideas concerning POD methods. The error analysis of the method is carried out in Section 5. Finally, some numerical experiments are shown in the last section to check in practice the theoretical results obtained.

\section{Preliminaries and notation}
As a model problem we consider the following reaction-diffusion equation 
\begin{equation}\label{eq:model}
\begin{array}{rcll}
u_t(t,x)-\nu\Delta u(t,x)+g(u(t,x))& = &f(t,x) &\quad (t,x)\in (0,T]\times \Omega,\\
u(t,x)&=&0,& \quad  (t,x)\in (0,T]\times \partial \Omega,\\
u(0,x)&=&u^0(x),&\quad  x\in \Omega,
\end{array}
\end{equation}
in a bounded domain $\Omega\subset {\Bbb R}^d$, $d\in \{2,3\}$, where $g$ is a nonlinear smooth function and $\nu$ is a positive parameter.
{For simplicity here and in the sequel we take $t_0=0$ as initial time}.

Let $C_p$ be the constant in the Poincar\'e inequality
\begin{equation}\label{poincare}
\|v\|_0\le C_p \|\nabla v\|_0,\quad v\in H_0^1(\Omega),
\end{equation}
{\color{red}where $\|\cdot\|_0$ denotes the standard $L^2$ norm}.
Let us denote by $X_h^k$ a finite element space based on piece-wise continuous polynomials of degree $k$  defined
over a triangulation of $\Omega$ of diameter $h$ 
and by $V_h^k$ the  finite element space based on piece-wise continuous polynomials of degree $k$ that satisfies also the homogeneous Dirichlet boundary conditions at the boundary $\partial \Omega$.

In the sequel,  $I_h u \in X_h^k$ will denote
the Lagrange interpolant of a continuous function $u$. The following bound can be found in \cite[Theorem 4.4.4]{brenner-scot}
\begin{equation}\label{cota_inter}
|u-I_h u|_{W^{m,p}(K)}\le c_{\rm int} h^{n-m}|u|_{W^{n,p}(K)},\quad 0\le m\le n\le k+1,
\end{equation}
where $n>d/p$ when $1< p\le \infty$ and $n\ge d$ when $p=1$.

Let us assume we want to approximate (\ref{eq:model}) in a given time interval $[0,T]$. We consider the semi-discrete finite element approximation: Find $u_h\ :\ (0,T]\rightarrow V_h^k$ such that
\begin{equation}\label{gal_semi}
(u_{h,t},v_h)+\nu(\nabla u_h,\nabla v_h)+(g(u_h),v_h)=(f,v_h),\quad \forall\ v_h\in V_h^k,
\end{equation}
with $u_h(0)=I_h u^0\in V_h^k.$ {\color{red}In (\ref{gal_semi}) $u_{h,t}$ denotes the derivative of $u_h$  with respect to time.}
If the weak solution of (\ref{eq:model}) is sufficiently smooth ({\color{red}$u, u_t\in H^{k+1}(\Omega)$}), then the following estimate {\color{red}for the error of the standard semi-discrete Galerkin approximation for $g$ locally Lipschitz continuous} is well known,
see for example \cite[Proof of Theorem 1]{bosco-titi-fem}, \cite[Theorem 14.1]{Thomee}, \cite[Lemma 4.2]{Wang}, 
\begin{equation}\label{cota_gal}
\max_{0\le t\le T}\left(\|(u-u_h)(t)\|_0+h\|(u-u_h)(t)\|_1\right)\le C(u)h^{k+1},
\end{equation}
{\color{red}where $\|\cdot\|_1$ denotes the standard norm in $H^1(\Omega)$.}
\section{BDF methods}
In this section we follow the notation and theory presented in \cite{bdf_stokes}. For the temporal integration of the reduced-order method introduced in Section 4 we are going to use the $q$-step backward differentiation formulae (BDF). For $T>0$, let
us denote by $\Delta t=T/M$ and by $t_n=n\Delta t$, $n=0,\ldots,M$ a uniform partition of $[0,T]$. For $u^n=u(t_n)$, the BDF approximation of order
$q$ to the time derivative $u_t$ is given by
\begin{equation}\label{eq_bdf_q}
\overline \partial_q u^n=\frac{1}{\Delta t}\sum_{i=0}^q \delta_i \color{red}{u^{n-i}},\ n\ge q,
\end{equation}  
{\color{red}where the method coefficients are determined from the relation (see \cite{Hairer})}
$$
 \delta(\zeta)=\sum_{i=0}^q \delta_i\zeta^i
=\sum_{l=1}^q\frac{1}{l}(1-\zeta)^l,
$$
We reproduce the results on $G$-stability for BDF schemes of \cite{bdf_stokes} that allow to use energy estimates in the
error analysis of the methods. The following lemma is stated and proved in \cite[Lemma 3.6]{bdf_stokes}, see also \cite{dalquis1}, \cite{dalquis2}.
\begin{lemma}\label{lema_3_6}
Let $\delta(\zeta)=\sum_{i=0}^q \delta_i \zeta^i$ and $\mu(\zeta)=\sum_{j=0}^q\mu_j\zeta^j$ be polynomials that have no common divisor. Assume the following condition holds
\begin{equation}\label{con_g}
{\rm Re}\frac{\delta(\zeta)}{\mu(\zeta)}>0,\quad \zeta\in {\Bbb C},|\zeta|<1.
\end{equation}
Let $(\cdot,\cdot)$
be a semi-inner product on a Hilbert space $H$ with associated norm $|\cdot|$. Then, there exists a symmetric positive-definite
matrix $G=[g_{i,j}]\in {\Bbb R}^{q\times q}$ such that for $v_0,\ldots,v_q\in H$, the following bound holds
\begin{equation}\label{g_stability}
{\rm Re}\left(\sum_{i=0}^q \delta_i v^{q-i},\sum_{j=0}^{\color{red} q}\mu_j v^{q-j}\right)\ge \sum_{i,j=1}^qg_{i,j}(v^i,v^j)-\sum_{i,j=1}^q g_{i,j}(v^{i-1},v^{j-1}).
\end{equation}
\end{lemma}
The application of $G$-stability to BDF schemes is ensured by the following lemma (see \cite[Lemma 3.7]{bdf_stokes}) based on the multiplier technique of Nevanlinna and Odeh \cite{nevan_ode}.
\begin{lemma} For $q=1,\ldots,5$ there exits $0\le \eta_q<1$ such that for $\delta(\zeta)=\sum_{l=1}^q (1/l)(1-\zeta)^l$ and $\mu(\zeta)=1-\eta_q\zeta$  condition (\ref{con_g}) holds. The values of $\eta_q$ ({\color{red}given up to four significant digits}) are found to be
\begin{equation}\label{eta}
\eta_1=\eta_2=0,\ \eta_3=0.0769,\ \eta_4=0.2878,\ \eta_5=0.8097.
\end{equation}
\end{lemma}

As in \cite{bdf_stokes} we define a $G$-norm associated with the semi-inner product $(\cdot,\cdot)$ on a Hilbert space $H$. Given
$\bV^n=[v^n,\ldots,v^{n-q+1}]$ with $v^{n-i+1}\in H$, $i=1,\ldots,q$, we define
\begin{equation}\label{norma_G}
\left|\bV^n\right|_G^2=\sum_{i,j=1}^q g_{i,j}(v^{n-i+1},v^{n-j+1}),
\end{equation}
where $G=[g_{i,j}]$ is the matrix in Lemma \ref{lema_3_6}. Since $G$ is symmetric positive definite the following inequality holds
\begin{equation}\label{equi_nor}
\tilde\lambda_m|v^n|^2\le \tilde\lambda_{\rm m}\sum_{j=1}^q|v^{n-j+1}|^2\le \left|\bV^n\right|_G^2\le \tilde\lambda_M\sum_{j=1}^q|v^{n-j+1}|^2,
\end{equation}
where $\tilde\lambda_m$ and $\tilde\lambda_M$ denotes the smallest and largest eigenvalues of $G$, respectively, and $|\cdot|$ is
the semi-norm induced by the semi-inner product $(\cdot,\cdot)$.

\section{Proper orthogonal decomposition}\label{sec:POD}

For $T>0$, $\Delta t=T/M$, $t_j=j\Delta t$, $j=0,\ldots,M$ and $N=M+1$, we define the 
space
\begin{equation}
\label{labU}
\bU = {\rm span}\left\{\sqrt{N}w_0,\tau\frac{u_h(t_1)-u_h(t_0)}{\Delta t},\tau\frac{u_h(t_2)-u_h(t_1)}{\Delta t},\ldots,\tau\frac{u_h(t_M)-u_h(t_{M-1})}{\Delta t}\right\},\\
\end{equation}
where $w_0$ is either $w_0=u_h(t_0)$ or $w_0=\overline u_h=\sum_{j=0}^Mu_h(t_j)/(M+1)$,
and $\tau$ is a time scale to make the snapshots dimensionally correct. {This means
that if $u_h$ is the finite element approximation to some physical quantity (velocity, concentration of a chemical reactant, pressure, etc) then, $w_0$ and $\tau({u_h(t_j)-u_h(t_{j-1})})/{\Delta t}$, $j=1,\ldots, M$, have the same physical units}. The
following analysis only requires $\tau>0$. 
We denote 
\begin{equation}\label{def_dif}
y_h^1=\sqrt{N}w_0, \qquad y_h^j = \tau\frac{u_h(t_{j-1})-u_h(t_{j-2})}{\Delta t}, \quad j=2,\ldots,M+1=N,
\end{equation}
so that $\bU={\rm span}\left\{y_h^1,y_h^2,\ldots,y_h^N\right\}$.

Let 
$X=H_0^1(\Omega)$ and let us denote the correlation matrix by $K=((k_{i,j}))\in {\Bbb R}^{N\times N}$ with
$$
k_{i,j}=\frac{1}{N}(y_h^i,y_h^j)_X=\frac{1}{N}(\nabla y_h^i,\nabla y_h^j),\quad i,j=1,\ldots,N,
$$
and $(\cdot,\cdot)$ the inner product in $L^2(\Omega)$. We denote by $\lambda_1\ge \lambda_2\ldots\ge \lambda_{d_r}>0$ the positive eigenvalues of $K$ and
by $\bv_1,\ldots,\bv_{d_r}\in {\Bbb R}^N$ the associated eigenvectors. The orthonormal POD basis functions of $\bU$ are
\begin{equation}\label{eq:varphi}
\varphi_k=\frac{1}{\sqrt{N}}\frac{1}{\sqrt{\lambda_k}}\sum_{j=1}^Nv_k^j y_h^j,
\end{equation}
where $v_k^j$ is the $j$ component of the eigenvector $\bv_k$.
For any $1\le r\le d_r$ let 
\begin{equation}\label{eq:bU_r}
\bU^r={\rm span}\left\{\varphi_1,\varphi_2,\ldots,\varphi_r\right\},
\end{equation}
and let us denote by $P^r:H_0^1(\Omega)\rightarrow \bU^r$ the $H_0^1$-orthogonal projection onto $\bU^r$. Then, it holds, {\color{red}see \cite[Proposition 1]{ku-Vol}},
\begin{equation}\label{cota_pod}
\frac{1}{N}\sum_{j=1}^N\|\nabla(y_h^j-P^r y_h^j)\|_0^2=\sum_{k={r+1}}^{d_r}\lambda_k,
\end{equation}
which,{ \textcolor{red}in view of (\ref{def_dif}}), gives
{\color{red}
\begin{eqnarray}\label{cota_pod2}
\|\nabla (I-P^r)w_0\|_0^2+\frac{1}{M+1}\frac{\tau^2}{\Delta t^2}\sum_{j=1}^M\|\nabla (I-P^r)D u_h(t_j)\|_0^2=\sum_{k={r+1}}^{d_r}\lambda_k,
\end{eqnarray}
}
for
$$
D u_h(t_j)=u_h(t_j)-u_h(t_{j-1}).
$$
{\color{red}The definition of the set of snapshots in (\ref{def_dif}) allows to get pointwise estimates.} The following lemma is taken from \cite{temporal_nos}.
\begin{lemma}\label{le:maxPr_FD} Let $\tilde C=1$ if $w_0=u_h(0)$ and $\tilde C=4$ if~$w_0=\overline u_h$.
{Let us denote by $u_h^n=u_h(t_n)$}. 
The following bound holds 
\begin{eqnarray}\label{max_dif0}
\max_{0\le n\le M}\|u_h^n-P^ru_h^n\|_0^2&\le& \left(2+4\tilde C \frac{T^2}{\tau^2}\right)C_p^2\sum_{k={r+1}}^{d_r}\lambda_k,\\
\label{max_dif1}
\max_{0\le n\le M}\|\nabla(u_h^n-P^ru_h^n)\|_0^2&\le& \left(2+4\tilde C \frac{T^2}{\tau^2}\right)\sum_{k={r+1}}^{d_r}\lambda_k.
\end{eqnarray}
As a consequence, since $\Delta t = T/M$ and  $(M+1)/M\le 2$, it follows that
\begin{equation}\label{max_dif_promedio}
{\Delta t }\sum_{n=0}^M \|u_h^n-P^ru_h^n\|_0^2\le T\left(4+8\tilde C \frac{T^2}{\tau^2}\right)C_p^2\sum_{k={r+1}}^{d_r}\lambda_k.
\end{equation}
\end{lemma}
\section{Error analysis of the POD method}
 We will consider the following fully-discrete POD-ROM approximation to approach (\ref{eq:model}): Find $u_r:(0,T]\rightarrow \bU^r$
such that
\begin{equation}\label{eq:pod}
(\overline \partial_q u_r^n,v_r)+\nu(\nabla u_r^n,\nabla v_r)+(g(u_r^n),v_r)=(f^n,v_r),\quad \forall\ v_r\in \bU^r, \quad n\ge q,
\end{equation}
with $u_r^0,\ldots,u_r^{q-1}\in \bU^r$.

For the error analysis of this section, we will assume that the nonlinear function
$g : \mathbb R \to \mathbb R$ is Lipschitz continuous (globally), with $g(u_r)$ being a short form  
notation of $g(u_r(t,x))$, $u_r(t,x)\in 
\mathbb R$, $t\in [0,T],\  x\in \Omega$. The extension of the error analysis to a more general function can be carried out arguing
as in \cite[Theorem 2]{temporal_nos}. 

Concerning the initial conditions $u_r^0,\ldots,u_r^{q-1}$ we will assume for simplicity that $u_r^j=P^r u_h(t_j)$, $0\le j\le q-1$.
In practice one could also take $u_r^0=P^r u_0$ and then obtain  $u_r^j$, $1\le j\le q-1,$ using as a time integrator in  the POD method 
a one step method of order $q$. There are also other options, see \cite[Remark 3.4]{bdf_stokes}.

We observe that to get the POD approximation (\ref{eq:pod}) one has to solve a non-linear system of equations. It is easy to show existence and uniqueness of solution for $\Delta t$ small enough. In practice, one can apply Newton's method to get an approximation
to the true solution.

For the error analysis we will apply the following discrete Gronwall inequality taken from \cite[Lemma 3.2]{nos_bdf}.
\begin{lemma}\label{gronwall2}
Let $B,a_j,b_j,c_j,\gamma_j,\color{red}{d_j}$ be nonnegative numbers such that
$$
a_n+\sum_{j=0}^n b_j\le \gamma_n a_n +\sum_{j=0}^{n-1} (\gamma_j+{\color{red}d_j}) a_j+\sum_{j=0}^n c_j+B,\quad for \quad n\ge 0.
$$
Suppose that $\gamma_j<1$, for all j, and set $\sigma_j=(1-\gamma_j)^{-1}$. Then
$$
a_n+\sum_{j=0}^nb_j\le \exp\left(\sigma_n \gamma_n+\sum_{j=0}^{n-1}(\sigma_j\gamma_j+{\color{red}d_j})\right)\left\{\sum_{j=0}^nc_j+B\right\}, \quad for\quad  n\ge 0.
$$
\end{lemma}
In the sequel, we will bound the error between the fully discrete POD approximation of (\ref{eq:pod}), $u_r^n$, and the projection
of the finite element approximation at time $t_n$, $P^r u_h^n$. 

Taking into account that we consider
projections respect to the $H_0^1$ inner product, that $\bU^r\subset V_h^k$ and using (\ref{gal_semi}), it is easy to see that $P^r u_h^n$ satisfies the following equation
\begin{eqnarray}\label{eq_pr}
(\overline \partial_q P^r u_h^n,v_r)+\nu(\nabla P^r u_h^n,\nabla v_r)+(g(u_h^n),v_r)=(f^n,v_r)+(\overline \partial_q P^r u_h^n-u_{h,t}^n,v_r), \forall v_r\in \bU^r.
\end{eqnarray}
\begin{theorem} \label{th1} ({\it Stability}) Let $u_h^n$ be the finite element approximation at time $t_n$ and $u_r^n$ be the POD approximation at time $t_n$.
Let us denote by $L$ the Lipchitz constant of $g$, assuming the time step $\Delta t$ satisfies
{\color{red}
\begin{eqnarray}\label{cond_deltat1}
&&\Delta t\left({1}+{5L}\right)\le {\delta_0},\quad \Delta t\frac{(3L+2{\eta_q})}{\tilde\lambda_m}\le 1,
\end{eqnarray}
}
with $\delta_0$ defined in (\ref{eq_bdf_q}),
then, there exists a constant $C>0$ such that the following bound holds
\begin{eqnarray}\label{er_defi_le}
\|u_r^M-P^r u_h^M\|_0^2+\sum_{n=q}^M\Delta t \nu \|\nabla (u_r^n-P^r u_h^n)\|_0^2\le C e^{2\alpha T} \left(\sum_{n=q}^M\Delta t \|\tau_r^n\|_0^2\right),
\end{eqnarray}
where $\color{red}{\alpha=L{(3+2{\eta_q})}/(2{\tilde\lambda_m})}$ and
\begin{eqnarray}\label{trun2}
&&\|\tau_r^q\|_0^2=c_q {\color{red}L}\|(I-P^r)u_h^q\|_0^2+c_q\|u_{h,t}^q-\overline \partial_q P^r u_h^q\|_0^2,\\
\label{trun1}
&&\|\tau_r^n\|_0^2=\frac{{\color{red}L}(1+2\eta_q)}{2}\|(I-P^r)u_h^n\|_0^2+\frac{(1+\eta_q)}{2}\|u_{h,t}^n-\overline \partial_q P^r u_h^n\|_0^2,\ n=q+1,\ldots,M,\qquad
\end{eqnarray}
and $c_q$ is a constant depending on $q$, $\delta_i$, $\tilde\lambda_M$ and $\eta_q$, specified along the proof of the theorem.
\end{theorem}
\begin{proof}
Let us denote by
$$
e_r^n=u_r^n-P^r u_h^n.
$$
Subtracting (\ref{eq_pr}) from (\ref{eq:pod}) we obtain the following error equation.
\begin{eqnarray}\label{eq_er}
(\overline \partial_q e_r^n,v_r)+\nu (\nabla e_r^n,\nabla v_r)=(g(u_h^n)-g(u_r^n),v_r)+(u_{h,t}^n-\overline \partial_q P^r u_h^n,v_r), \ \forall v_r\in \bU^r.
\end{eqnarray}
Now, we take in (\ref{eq_er})
$v_r=e_r^n-\eta_q e_r^{n-1}$ to obtain
\begin{eqnarray}\label{eq_er2}
(\overline \partial_q e_r^n,e_r^n-\eta_q e_r^{n-1})+\nu \|\nabla e_r^n\|_0^2&=&\eta_q
\nu (\nabla e_r^n,\nabla e_r^{n-1})+(g(u_h^n)-g(u_r^n),e_r^n-\eta_q e_r^{n-1})\nonumber\\
&&\quad +(u_{h,t}^n-\overline \partial_q P^r u_h^n,e_r^n-\eta_q e_r^{n-1}).
\end{eqnarray}
Applying now Lemma \ref{lema_3_6} with $H=L^2(\Omega)$ and definition (\ref{norma_G}), from (\ref{eq_er2}) we get
\begin{eqnarray}\label{eq_er3}
&&\color{red}{|E_r^n|_G^2-|E_r^{n-1}|_G^2+\Delta t \nu \|\nabla e_r^n\|_0^2 \le \eta_q
\Delta t \nu (\nabla e_r^n,\nabla e_r^{n-1})}\nonumber\\
&&\quad +\Delta t (g(u_h^n)-g(u_r^n),e_r^n-\eta_q e_r^{n-1})+\Delta t (u_{h,t}^n-\overline \partial_q P^r u_h^n,e_r^n-\eta_q e_r^{n-1}).
\end{eqnarray}
For the first term on the right hand side of (\ref{eq_er3}) we write
\begin{eqnarray}\label{au_er1}
\eta_q\nu (\nabla e_r^n,\nabla e_r^{n-1})\le \frac{\eta_q}{2}\left(\nu\|\nabla e_r^n\|_0^2+\nu\|\nabla e_r^{n-1}\|_0^2\right)
\end{eqnarray}
For the next term we observe that, using that $g$ is Lipschitz continuous we get
$$
\|g(u_h^n)-g(u_r^n)\|_0\le L\|u_h^n-u_r^n\|_0\le L \left(\|(I-P^r)u_h^n\|_0+\|e_r^n\|_0\right).
$$
Then, using Young's inequality, we obtain
\begin{eqnarray}\label{au_er2}
&&(g(u_h^n)-g(u_r^n),e_r^n-\eta_q e_r^{n-1})\le \|g(u_h^n)-g(u_r^n)\|_0\|e_r^n\|_0
+\|g(u_h^n)-g(u_r^n)\|_0\eta_q\|e_r^{n-1}\|_0\nonumber\\
&&\quad \le L(\|e_r^n\|_0^2+\eta_q\|e_r^n\|_0\|e_r^{n-1}\|_0)+L\|(I-P^r)u_h^n\|_0(\|e_r^n\|_0+\eta_q\|e_r^{n-1}\|_0)\nonumber\\ &&\quad \le \color{red}{L\|e_r^n\|_0^2+L\eta_q\|e_r^n\|_0^2+\frac{L\eta_q}{4}\|e_r^{n-1}\|_0^2
+\frac{L}{2}\|(I-P^r)u_h^n\|_0^2+\frac{L}{2}\|e_r^n\|_0^2}\nonumber\\
&&\quad\quad+\color{red}{{L}\eta_q\|(I-P^r)u_h^n\|_0^2+\frac{L\eta_q}{4}\|e_r^{n-1}\|_0^2}\nonumber\\
&&\quad \le \color{red}{L\left(\frac{3}{2}+{\eta_q}\right)\|e_r^n\|_0^2+\frac{L\eta_q}{2}\|e_r^{n-1}\|_0^2+\frac{L(1+2\eta_q)}{2}\|(I-P^r)u_h^n\|_0^2}.
\end{eqnarray}
and
\begin{eqnarray}\label{au_er3}
(u_{h,t}^n-\overline \partial_q P^r u_h^n,e_r^n-\eta_q e_r^{n-1})&\le& \|u_{h,t}^n-\overline \partial_q P^r u_h^n \|_0\|e_r^n\|_0+\|u_{h,t}^n-\overline \partial_q P^r u_h^n\|_0\eta_q\|e_r^{n-1}\|_0\nonumber\\
&\le&\frac{(1+\eta_q)}{2}\|u_{h,t}^n-\overline \partial_q P^r u_h^n\|_0^2+\frac{1}{2}\|e_r^n\|_0^2+\frac{\eta_q}{2}\|e_r^{n-1}\|_0^2.
\end{eqnarray}
Inserting (\ref{au_er1}), (\ref{au_er2}) and (\ref{au_er3}) into (\ref{eq_er3}) we get
\begin{eqnarray}\label{eq_er4}
&&|E_r^n|_G^2-|E_r^{n-1}|_G^2+\frac{2-\eta_q}{2}\Delta t \nu \|\nabla e_r^n\|_0^2-\frac{\eta_q}{2}\Delta t\nu \|\nabla e_r^{n-1}\|_0^2\le \Delta t \color{red}{L\left(\frac{3}{2}+{\eta_q}\right)\|e_r^n\|_0^2}\qquad\nonumber\\
&&\quad +\Delta t {\color{red}L}\eta_q\|e_r^{n-1}\|_0^2+\Delta t \frac{\color{red}{L}(1+2\eta_q)}{2}\|(I-P^r)u_h^n\|_0^2+\Delta t\frac{(1+\eta_q)}{2}\|u_{h,t}^n-\overline \partial_q P^r u_h^n\|_0^2.
\end{eqnarray}
Now, adding terms in (\ref{eq_er4}) from $n=q+1,\ldots,M$ and using definition (\ref{trun1}) we obtain
\begin{eqnarray}\label{eq_er4}
&&|E_r^M|_G^2+(1-\eta_q)\sum_{n=q+1}^M\Delta t \nu \|\nabla e_r^n\|_0^2\le |E_r^q|_G^2+\frac{\eta_q}{2}\Delta t\nu \|\nabla e_r^{q}\|_0^2+\Delta t \color{red}{L}\eta_q\|e_r^{q}\|_0^2\nonumber\\
&&\quad +\Delta t {\color{red}L\left(\frac{3}{2}+\eta_q\right)\|e_r^M\|_0^2}+\sum_{n=q+1}^{M-1}\Delta t {\color{red}L\left(\frac{3}{2}+2{\eta_q}\right)}\|e_r^n\|_0^2+\sum_{n=q+1}^M\Delta t \|\tau_r^n\|_0^2.
\end{eqnarray}
Let us observe that $1-\eta_q>0$.

Applying (\ref{equi_nor}) to (\ref{eq_er4}) we obtain
\begin{eqnarray}\label{eq_er42}
&&\tilde\lambda_m\|e_r^M\|_0^2+(1-\eta_q)\sum_{n=q+1}^M\Delta t \nu \|\nabla e_r^n\|_0^2\le \tilde\lambda_M\sum_{n=1}^q\|e_r^n\|_0^2+\frac{\eta_q}{2}\Delta t\nu \|\nabla e_r^{q}\|_0^2+\Delta t {\color{red}L}\eta_q\|e_r^{q}\|_0^2\nonumber\\
&&\quad +\Delta t {\color{red}L}\left(\frac{3}{2}+{\eta_q}\right)\|e_r^M\|_0^2+\sum_{n=q+1}^{M-1}\Delta t {\color{red}L\left(\frac{3}{2}+2{\eta_q}\right)\|e_r^n\|_0^2}
+\sum_{n=q+1}^M\Delta t \|\tau_r^n\|_0^2.
\end{eqnarray}
For the first 3 terms on the right-hand side of (\ref{eq_er42}) we can write
\begin{eqnarray}\label{auxi}
\tilde\lambda_M\sum_{n=1}^q\|e_r^n\|_0^2+\frac{\eta_q}{2}\Delta t\nu \|\nabla e_r^{q}\|_0^2+\Delta t \eta_q\|e_r^{q}\|_0^2
&\le& c_{q,1}\|e_r^{q}\|_0^2+c_{q,2}\Delta t\nu \|\nabla e_r^{q}\|_0^2\nonumber\\
&&+\tilde\lambda_M\sum_{n=1}^{q-1}\|e_r^n\|_0^2,
\end{eqnarray}
where the constants $c_{q,1}$ and $c_{q,2}$ are defined as
$$
c_{q,1}=\tilde\lambda_M+\Delta t {\color{red}L}\eta_q, \ c_{q,2}=\frac{\eta_q}{2}.
$$
We need to bound the first two-terms on the right-hand side of (\ref{auxi}). To this end, we consider (\ref{eq_er}) for $n=q$ and we take $v_r=e_r^q$. Then, we get
\begin{eqnarray}\label{er_q1}
\delta_0 \|e_r^q\|_0^2+\Delta t \nu \|\nabla e_r^q\|_0^2&\le& \Delta t\|g(u_h^q)-g(u_r^q)\|_0\|e_r^q\|_0+\Delta t\|u_{h,t}^q-\overline \partial_q P^r u_h^q\|_0\|e_r^q\|_0\nonumber\\
&&+\left\|\sum_{i=1}^{q}\delta_ie_r^{q-i}\right\|_0\|e_r^q\|_0.
\end{eqnarray}
For the first term on the right-hand side of (\ref{er_q1}), arguing as before, we get
\begin{eqnarray}\label{er_q2}
\Delta t\|g(u_h^q)-g(u_r^q)\|_0\|e_r^q\|_0&\le& \Delta t L(\|(I-P^r)u_h^q\|_0+\|e_r^q\|_0)\|e_r^q\|_0\nonumber\\
&\le& \Delta t{\color{red}L}\|(I-P^r)u_h^q\|_0^2+{\color{red}\frac{5}{4}L\Delta t}\|e_r^q\|_0^2.
\end{eqnarray}
For the second term, Young's inequality yields
\begin{eqnarray}\label{er_q3}
\Delta t\|u_{h,t}^q-\overline \partial_q P^r u_h^q\|_0\|e_r^q\|_0\le \Delta t \|u_{h,t}^q-\overline \partial_q P^r u_h^q\|_0^2
+\frac{\Delta t}{4}\|e_r^q\|_0^2.
\end{eqnarray}
Finally, for the last term, applying again Young's inequality, we write
\begin{eqnarray}\label{er_q4}
\left\|\sum_{i=1}^{q}\delta_ie_r^{q-i}\right\|_0\|e_r^q\|_0\le \frac{1}{\delta_0}\left\|\sum_{i=1}^{q}\delta_ie_r^{q-i}\right\|
_0^2+\frac{\delta_0}{4}\|e_r^q\|_0^2\le c_{q,3}\sum_{i=0}^{q-1}\|e_r^i\|_0^2+\frac{\delta_0}{4}\|e_r^q\|_0^2,
\end{eqnarray}
where $c_{q,3}$ is a constant that depends on $q$ and $\delta_i$, $i=0,\ldots,q$. Now, we assume $\Delta t$ is small enough 
so that the first condition in (\ref{cond_deltat1}) holds. Inserting (\ref{er_q2}), (\ref{er_q3}) and (\ref{er_q4}) into (\ref{er_q1}) yields
\begin{eqnarray*}\label{er_q5}
\frac{\delta_0}{2} \|e_r^q\|_0^2+\Delta t \nu \|\nabla e_r^q\|_0^2&\le&\Delta t\left({\color{red}L}\|(I-P^r)u_h^q\|_0^2+\|u_{h,t}^q-\overline \partial_q P^r u_h^q\|_0^2\right)+c_{q,3}\sum_{i=0}^{q-1}\|e_r^i\|_0^2.\quad
\end{eqnarray*}
Then, from the above inequality, going back to (\ref{auxi}) it is easy to see that
\begin{eqnarray}\label{auxi2}
c_{q,1}\|e_r^{q}\|_0^2+c_{q,2}\Delta t\nu \|\nabla e_r^{q}\|_0^2\le c_q\Delta t \left({\color{red}L}\|(I-P^r)u_h^q\|_0^2+\|u_{h,t}^q-\overline \partial_q P^r u_h^q\|_0^2\right)+c_q \sum_{i=0}^{q-1}\|e_r^i\|_0^2,
\end{eqnarray}
where $c_q$ is a constant that depends on $\delta_0$, $c_{q,i}$, $i=1,2,3$.

Using definition (\ref{trun2}) and inserting (\ref{auxi}) and (\ref{auxi2}) into (\ref{eq_er42}) we reach
\begin{eqnarray}\label{eq_erd}
&&\|e_r^M\|_0^2+\frac{(1-\eta_q)}{\tilde\lambda_m}\sum_{n=q+1}^M\Delta t \nu \|\nabla e_r^n\|_0^2\le \frac{(c_q+\tilde\lambda_M)}{\tilde\lambda_m} \sum_{i=0}^{q-1}\|e_r^i\|_0\nonumber\\
&&\quad +\Delta t {\color{red}L\frac{(3+{2\eta_q})}{2\tilde\lambda_m}}\|e_r^M\|_0^2+\sum_{n=q+1}^{M-1}\Delta t {\color{red}L\frac{(3+4{\eta_q})}{2\tilde\lambda_m}}\|e_r^n\|_0^2
+\sum_{n=q}^M\frac{\Delta t}{\tilde\lambda_m} \|\tau_r^n\|_0^2.
\end{eqnarray}
Now, we can apply Lemma~\ref{gronwall2} with $a_n=\|e_r^n\|_0^2$, $\color{red}{\gamma_n=\Delta tL{(3+2{\eta_q})}/(2{\tilde\lambda_m})}$, $n=q+1,\ldots,M$, $\color{red}{d_n}={\color{red}L}\Delta t\eta_q/\tilde\lambda_m$, $n=q+1,\ldots,M-1$. Assuming the second condition in (\ref{cond_deltat1}) holds,
then $\gamma_n\le 1/2$
and $\sigma_n=(1-\gamma_n)^{-1}\le 2$. Then,
\begin{eqnarray}\label{eq_def}
&&\|e_r^M\|_0^2+\frac{(1-\eta_q)}{\tilde\lambda_m}\sum_{n=q+1}^M\Delta t \nu \|\nabla e_r^n\|_0^2 \nonumber\\
&&\ \le\exp\left(\sigma_M\gamma_M+\sum_{n=q+1}^{M-1}(\sigma_j\gamma_j+\delta_j)\right)\left(\frac{(c_q+\tilde\lambda_M)}{\tilde\lambda_m} \sum_{i=0}^{q-1}\|e_r^i\|_0+\sum_{n=q}^M\frac{\Delta t}{\tilde\lambda_m} \|\tau_r^n\|_0^2\right)\nonumber\\
&&\ \le\exp(\beta_M\Delta t)\left(\frac{(c_q+\tilde\lambda_M)}{\tilde\lambda_m} \sum_{i=0}^{q-1}\|e_r^i\|_0+\sum_{n=q}^M\frac{\Delta t}{\tilde\lambda_m} \|\tau_r^n\|_0^2\right),
\end{eqnarray}
where $\beta_M=2\alpha+(M-1-q)(2\alpha+\color{red}{L}\eta_q/\tilde\lambda_m)$, for $\color{red}{\alpha=L{(3+2{\eta_q})}/(2{\tilde\lambda_m})}$. We can finally simplify
(\ref{eq_def}) and write (taking into account (\ref{er_q5}) to bound $\|\nabla e_r^q\|_0^2$)
\begin{eqnarray}\label{er_defi}
\|e_r^M\|_0^2+\sum_{n=q}^M\Delta t \nu \|\nabla e_r^n\|_0^2\le C e^{2\alpha T} \left(\sum_{i=0}^{q-1}\|e_r^i\|_0+\sum_{n=q}^M\Delta t \|\tau_r^n\|_0^2\right),
\end{eqnarray}
for a constant $C$ depending on previous constants ($\tilde\lambda_m$, $\tilde\lambda_M$, $\eta_q$, $c_q$, \ldots). Finally, taking into account that $u_r^i=P^r u_h^i$, $i=0,\ldots,q-1$ we reach (\ref{er_defi_le}).
\end{proof}
\begin{theorem} \label{th2} ({\it Consistency}) Let $\tau_r^n$, $n=q,\ldots,M$ be the truncation errors defined in (\ref{trun2})-(\ref{trun1}).
Then, the following bound holds
\begin{eqnarray}\label{eq_trun_def}
\sum_{n=q}^M\Delta t \|\tau_r^n\|_0^2&\le& T\left(4+8\tilde C \frac{T^2}{\tau^2}+\frac{C}{\tau^2}\right)C_p^2\sum_{k={r+1}}^{d_r}\lambda_k+C (\Delta t)^{2q}\int_0^T\|\partial_t^{q+1}u_h(s)\|_0^2 \ ds,\quad
\end{eqnarray}
where $\tilde C=1$ if $w_0=u_h(0)$ and $\tilde C=4$ if~$w_0=\overline u_h$ and $C$ is a positive constant.
\end{theorem}
\begin{proof}
In view of definitions (\ref{trun2})-(\ref{trun1}) we only need to bound
\begin{eqnarray}\label{cota_trun}
\sum_{n=q}^M\Delta t \|(I-P^r)u_h^n\|_0^2,\quad \sum_{n=q}^M\Delta t\|u_{h,t}^n-\overline \partial_q P^r u_h^n\|_0^2.
\end{eqnarray}
For the first term on the right-hand side of (\ref{cota_trun}) we apply (\ref{max_dif_promedio}) to bound it in terms of the tail of the eigenvalues.
\begin{equation}\label{nece}
\sum_{n=q}^M {\Delta t }\|u_h^n-P^ru_h^n\|_0^2\le \sum_{n=0}^M {\Delta t }\|u_h^n-P^ru_h^n\|_0^2\le T\left(4+8\tilde C \frac{T^2}{\tau^2}\right)C_p^2\sum_{k={r+1}}^{d_r}\lambda_k.
\end{equation}
To bound the second term on the right-hand side of (\ref{cota_trun}) we decompose
\begin{eqnarray}\label{tiempo_1}
\|u_{h,t}^n-\overline \partial_q P^r u_h^n\|_0^2\le 2\|u_{h,t}^n-\overline \partial_q  u_h^n\|_0^2
+2\|(I-P^r)\overline \partial_q  u_h^n\|_0^2.
\end{eqnarray}
To bound the first term on the right-hand side above we apply \cite[(B.2)]{bdf_stokes}
\begin{eqnarray*}
\|u_{h,t}^n-\overline \partial_q  u_h^n\|_0\le C(\Delta t)^{q-1/2}\left(\int_{t_{n-q}}^{t_n}\|\partial_t^{q+1}u_h(s)\|_0^2 \ ds\right)^{1/2},
\end{eqnarray*}
for a positive constant $C$.
Then, we get
\begin{eqnarray}\label{tiempo2}
\sum_{n=q}^M\Delta t \|u_{h,t}^n-\overline \partial_q  u_h^n\|_0^2\le C (\Delta t)^{2q}\int_0^T\|\partial_t^{q+1}u_h(s)\|_0^2 \ ds.
\end{eqnarray}

To conclude, we need to bound the second term on the right-hand side of (\ref{tiempo_1}). The key point to get this bound is to observe that $\overline \partial_q  u_h^n$ can be written as a linear combination of first order finite differences as those
appearing in (\ref{cota_pod2}).

From (\ref{eq_bdf_q}) we can write
$$
\overline \partial_q  u_h^n=\frac{1}{\Delta t}\left(\delta_0 u_h^{n}+\delta_{1}u_h^{n-1}+\ldots+\delta_{q-1}u_h^{n-q+1}+\delta_q u_h^{n-q}\right).
$$
Now, we observe that
\begin{eqnarray*}
&&\delta_0 u_h^{n}+\delta_{1}u_h^{n-1}+\ldots+\delta_{q-1}u_h^{n-q+1}+\delta_q u_h^{n-q}=\nonumber\\
&&\delta_0 Du_h^{n}+(\delta_0+\delta_{1})u_h^{n-1}+\ldots+\delta_{q-1}u_h^{n-q+1}+\delta_q u_h^{n-q}=\nonumber\\
&&\delta_0 Du_h^{n}+(\delta_0+\delta_{1})D u_h^{n-1}+(\delta_0+\delta_{1}+\delta_2)u_h^{n-2}\ldots+\delta_{q-1}u_h^{n-q+1}+\delta_q u_h^{n-q}=\nonumber\\
&&\ldots\nonumber\\
&&\delta_0 D u_h^{n}+(\delta_0+\delta_{1})D u_h^{n-1}+\ldots+(\delta_0+\delta_{1}+\ldots+\delta_{q-1})u_h^{n-q+1}+\delta_q u_h^{n-q}.
\end{eqnarray*}
Taking into account that $\delta_0+\delta_1+\ldots+\delta_q=0$ by consistency (the BDF-q formula applied to the constant function 1
is an approximation to its derivative,  i.e., 0), we obtain
\begin{eqnarray*}
&&\delta_0 u_h^{n}+\delta_{1}u_h^{n-1}+\ldots+\delta_{q-1}u_h^{n-q+1}+\delta_q u_h^{n-q}=\nonumber\\
&&\delta_0 D u_h^{n}+(\delta_0+\delta_{1})D u_h^{n-1}+\ldots+(-\delta_q)u_h^{n-q+1}+\delta_q u_h^{n-q}\nonumber \\
&&\delta_0 D u_h^{n}+(\delta_0+\delta_{1})D u_h^{n-1}+\ldots-\delta_q D u_h^{n-q+1}\nonumber\\
&&\alpha_0 D u_h^{n}+\alpha_{1}D u_h^{n-1}+\ldots+\alpha_{q-1} D u_h^{n-q+1},
\end{eqnarray*}
where
$$
\alpha_0=\delta_0,\quad {\color{red}\alpha_{j}=\sum_{i=0}^j\delta_i},\ i=1,\ldots,q-2,\quad \alpha_{q-1}=-\delta_q,
$$
so that
\begin{equation}\label{key}
\overline \partial_q  u_h^n=\frac{1}{\Delta t}\sum_{j=0}^{q-1} \alpha_j D u_h^{n-j}.
\end{equation}
Using (\ref{key}) and denoting by $D^1=D/\Delta t $, we can write
\begin{eqnarray}\label{tiempo3}
&&\sum_{n=q}^M\Delta t\|(I-P^r)\overline \partial_q  u_h^n\|_0^2=\sum_{n=q}^M\Delta t\left\|(I-P^r)\sum_{j=0}^{q-1} \alpha_j D^1 u_h^{n-j}\right\|_0^2=\nonumber\\ 
&&\quad \sum_{n=q}^M\Delta t\left\|\sum_{j=0}^{q-1} \alpha_j (I-P^r) D^1 u_h^{n-j}\right\|_0^2\le C \sum_{n=1}^M\Delta t\|(I-P^r)D^1u_h^n\|_0^2,
\end{eqnarray}
for a positive constant $C$. 

Now, we observe that from (\ref{cota_pod2}) we obtain
$$
\frac{\tau^2}{M+1}\sum_{j=1}^M\|\nabla D^1 u_h(t_j)\|_0^2\le\sum_{k={r+1}}^{d_r}\lambda_k,
$$
And then, the above inequality, Poincar\'e inequality (\ref{poincare})
and (\ref{tiempo3})  yield
\begin{eqnarray}\label{tiempo4}
\sum_{n=q}^M\Delta t\|(I-P^r)\overline \partial_q  u_h^n\|_0^2\le C \frac{T}{\tau^2}C_p^2\sum_{k={r+1}}^{d_r}\lambda_k.
\end{eqnarray}
From (\ref{nece}), (\ref{tiempo_1}), (\ref{tiempo2}) and (\ref{tiempo4}) we conclude (\ref{eq_trun_def}).
\end{proof}
From Theorems \ref{th1} and \ref{th2} we finally obtain
\begin{theorem} \label{th3} ({\it Convergence}) Let $u_h^n$ be the finite element approximation at time $t_n$ and $u_r^n$ be the POD approximation at time $t_n$. Assuming the time step $\Delta t$ satisfies (\ref{cond_deltat1}),
then the following bound holds
\begin{eqnarray}\label{er_defi_le2}
&&\|u_r^M-P^r u_h^M\|_0^2+\sum_{n=q}^M\Delta t \nu \|\nabla (u_r^n-P^r u_h^n)\|_0^2\nonumber\\
&&\quad\le C e^{2\alpha T} \left(T\left(4+8\tilde C \frac{T^2}{\tau^2}+\frac{C}{\tau^2}\right)C_p^2\sum_{k={r+1}}^{d_r}\lambda_k+C (\Delta t)^{2q}\int_0^T\|\partial_t^{q+1}u_h(s)\|_0^2 \ ds\right),\quad
\end{eqnarray}
where $\color{red}{\alpha=L{(3+2{\eta_q})}/(2{\tilde\lambda_m})}$, $\tilde C=1$ if $w_0=u_h(0)$ and $\tilde C=4$ if~$w_0=\overline u_h$ and $C$ is a positive constant.
\end{theorem}
To conclude, we bound the error between the POD approximation and the true solution. We can prove pointwise-in-time error
estimations thanks to Lemma \ref{le:maxPr_FD} in which the use of difference quotients in the set of snapshots is essential.
\begin{theorem} \label{th4} ({\it Convergence}) Let $u$ be the solution of (\ref{eq:model}), $u_h$ the continuous-in-time finite element approximation and $u_r^n$ be the POD approximation at time $t_n$. Assuming the time step $\Delta t$ satisfies (\ref{cond_deltat1}),
then the following bound holds
\begin{eqnarray}\label{er_defi_le3}
&&\max_{q\le n\le M}\|u_r^n-u(t_n)\|_0^2+\sum_{n=q}^M\Delta t \nu \|\nabla (u_r^n-u(t_n))\|_0^2\nonumber\\
&&\quad\le C e^{2\alpha T} \left (T\left(1+\frac{T^2}{\tau^2}\right)C_p^2\sum_{k={r+1}}^{d_r}\lambda_k+(\Delta t)^{2q}\int_0^T\|\partial_t^{q+1}u_h(s)\|_0^2 \ ds\right)\nonumber\\
&&\quad+C\left(1+\frac{T^2}{\tau^2}\right)(C_p^2+T\nu)\sum_{k={r+1}}^{d_r}\lambda_k+C(u)(h^{k+1}+ T\nu h^k).
\end{eqnarray}
where $\color{red}{\alpha=L{(3+2{\eta_q})}/(2{\tilde\lambda_m})}$ and $C$ is a positive constant.
\end{theorem}
\begin{proof}
We will use the following decomposition of the error
\begin{equation}\label{decom}
u_r^n-u(t_n)=(u_r^n-P^r u_h^n)+(P^r u_h^n-u_h^n)+(u_h^n-u(t_n)).
\end{equation}
To bound the error of the first term on the right-hand side above we apply Theorem \ref{th3}. Let us observe that the
bound of $\|u_r^n-P^r u_h^n\|_0^2$ for $n=M$ in (\ref{er_defi_le2}) also holds for any other value of $n$ with $q\le n\le M$
since taking any $n$, $q\le n\le M$, in the proof of (\ref{er_defi_le2}) instead of $M$ does not change the proof.
Let us also observe that the error contribution of this term is included on the right-hand side of the error bound (\ref{er_defi_le3}).

To bound the error in the second term in (\ref{decom}) we apply (\ref{max_dif0}), (\ref{max_dif1}) from Lemma~\ref{le:maxPr_FD}, then
\begin{eqnarray*}
\max_{q\le n\le M}\|u_h^n-P^r u_h^n\|_0^2+\sum_{n=q}^M\Delta t \nu \|\nabla (u_h^n-P^r u_h^n)\|_0^2&\le&
\left(2+4\tilde C \frac{T^2}{\tau^2}\right)C_p^2\sum_{k={r+1}}^{d_r}\lambda_k\nonumber\\
&&\ +T\nu \left(2+4\tilde C \frac{T^2}{\tau^2}\right)\sum_{k={r+1}}^{d_r}\lambda_k.
\end{eqnarray*}
We observe that the right-hand side above is also included on the right-hand side of the error bound (\ref{er_defi_le3}).

Finally, the bound of the last term in (\ref{decom}) is easily obtained from (\ref{cota_gal}).
\begin{eqnarray}\label{snapshots}
\max_{q\le n\le M}\|u_h^n-u(t_n)\|_0^2+\sum_{n=q}^M\Delta t \nu \|\nabla (u_h^n-u(t_n))\|_0^2\le C(u)(h^{k+1}+ T\nu h^k).
\end{eqnarray}
\end{proof}
{\color{red}
\begin{remark}
Let us observe that the  last term in the bound (\ref{er_defi_le3}) accounts for the error in the computation of the snapshots, see (\ref{snapshots}).  In this paper we are assuming that the snapshots $u_h(t_n)$ do not have temporal error since we have assumed that $u_h(t)$ is the continuous in time semi-discrete Galerkin approximation. In practice, there is always a temporal error in the computation of the snapshots that would appear added to the right-hand side of (\ref{snapshots}), and, as a consequence, also to the right-hand side of (\ref{er_defi_le3}). In the practical implementation of the method, the snapshots are computed offline and we assume that they are sufficiently accurate for the end user's requirements. 
\end{remark}
}

\begin{remark} The final bound (\ref{er_defi_le3}) for the error is optimal in terms of the tail of the eigenvalues (\ref{the_sigma_r}), the temporal
error and the spatial discretization error. 

We observe that, assuming that the true solution is smooth enough in time (which is needed to get rate of convergence $q$ in time)
 the quantity $\int_0^T\|\partial_t^{q+1}u_h(s)\|_0^2 \ ds$ can be bounded as
$$
\int_0^T\|\partial_t^{q+1}u_h(s)\|_0^2 \ ds \le 2 \int_0^T\|\partial_t^{q+1}(u_h-u)(s)\|_0^2+\int_0^T\|\partial_t^{q+1}u(s)\|_0^2 \ ds,
$$
The first term can be bounded taking time derivatives respect to $t$ in (\ref{eq:model}) and (\ref{gal_semi}) and then
applying the same kind of techniques that give (\ref{cota_gal}).
\end{remark}
\bigskip

\begin{remark} The error analysis corresponding to the BDF 6 method can also be obtained applying the ideas in \cite{bdf_stokes} to the error analysis shown in the present paper. The analysis of the time discretization error of BDF 6 is however more complicated and has to be done separately from the analysis corresponding to $q\le 5$. Since the key ideas for dealing with the time discretization error can already be found in
\cite{bdf_stokes} we think that there is no need to include this analysis here.
\end{remark}
\bigskip

\begin{remark}
As stated in the introduction, the error analysis of the case in which the standard set of snapshots based on the computed finite element approximations at different times are used, instead of difference quotients, can be carried out following \cite{Bosco_pointwise}. In that case, the error bounds would not be optimal in terms of the tail of the eigenvalues (\ref{the_sigma_r}), but as close to 
optimal as the smoothness of the function providing the snapshots allow. The analogous bound to (\ref{cota_pod2}) for 
the error projections onto the POD space of the difference quotients is also shown in \cite{Bosco_pointwise}, being again
as close to 
optimal as the smoothness of the function providing the snapshots allow.
\end{remark}
\section{Numerical Experiments}

We consider the system known as the Brusselator with diffusion
\begin{equation} \label{bruss}
	\begin{array}{rclcl}
		u_t&=&\nu\Delta u +1 +u^2v  - 4 u,&\qquad& (x,t)\in \Omega \times (0,T],\\
		v_t &=&\nu \Delta v+3 u -u^2v,& \qquad
		& (x,t)\in \Omega \times (0,T],\\
		&&u(x,t)=1, \quad v(x,t) =3,&& (x,t)\in \Gamma_1\times(0,T],\\
		&&\partial_n u(x,t)=\partial_n v(x,t) =0,&&  (x,t)\in \Gamma_2\times(0,T],\\
	\end{array}
\end{equation}
where $\nu$ is a positive parameter, $\Omega=[0,1]\times[0,1]$, $\Gamma_1\subset\partial\Omega$ is the union of of sides
$\{x=1\}\bigcup \{ y=1\}$, and~$\Gamma_2$ is the rest of~$\partial\Omega$. This system has an unstable equilibrium $u=1$, $v=3$, and, for $\nu$, sufficiently small, a stable limit cycle (see e.g., \cite[\S~I.16,\S~2.10]{Hairer}). If the boundary conditions are $\partial_n u(x,t)=\partial_n v(x,t) =0$ for $(x,t)\in(0,T]\times\partial\Omega$, then, the periodic orbit is flat in space, $\nabla u=\nabla v=0$, but with those in~(\ref{bruss}) it has more spatial complexity. This system was used in the numerical experiments of \cite[Section~6]{temporal_nos}. 

\begin{figure}[htbp]
	\centering
	\begin{minipage}{\textwidth}
		\centering
		\includegraphics[width=0.24\linewidth]{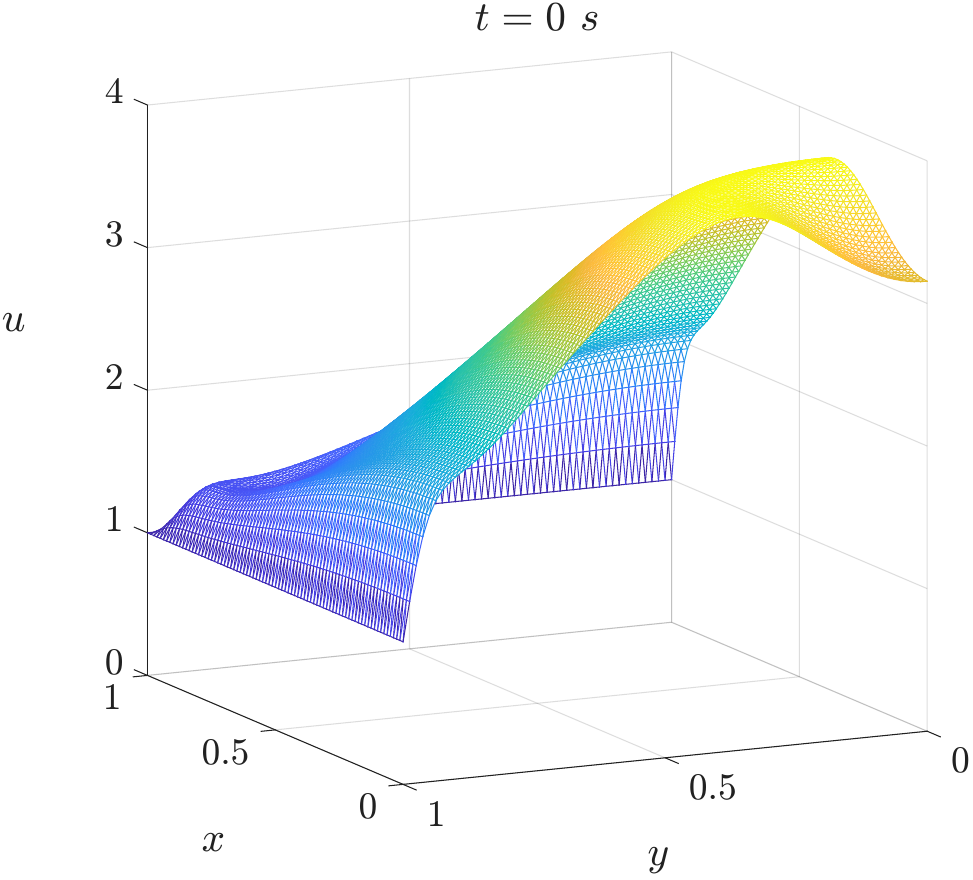}
		\includegraphics[width=0.24\linewidth]{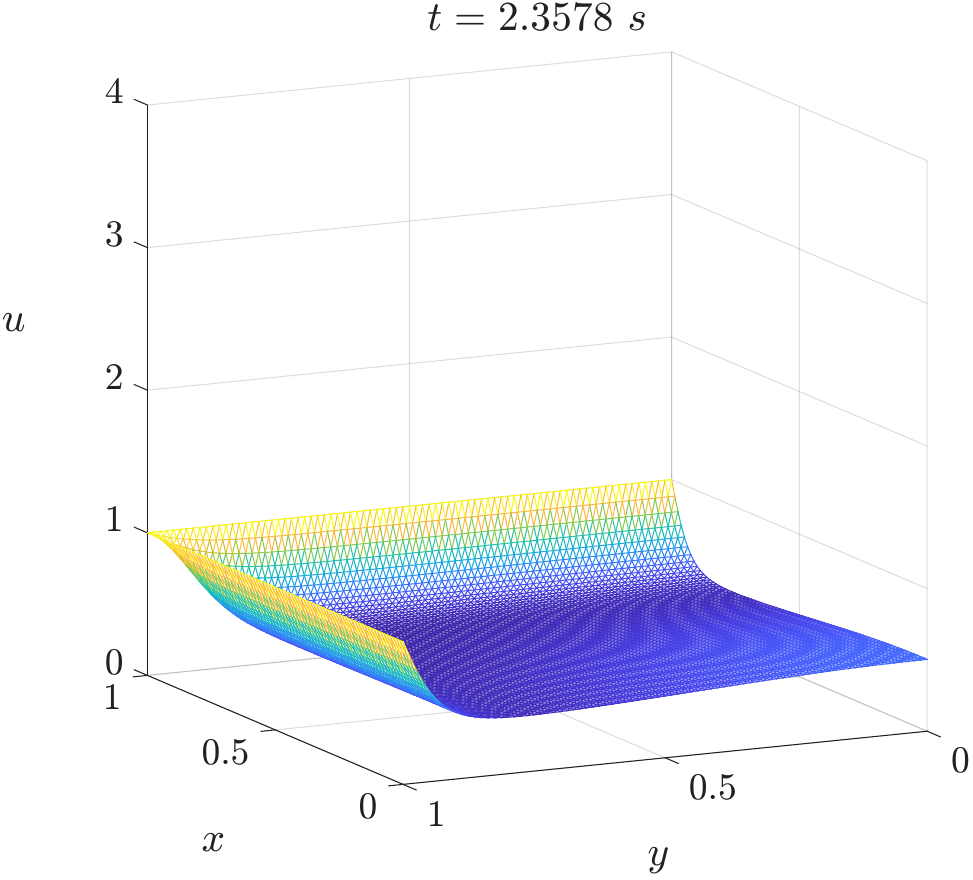} 
		\includegraphics[width=0.24\linewidth]{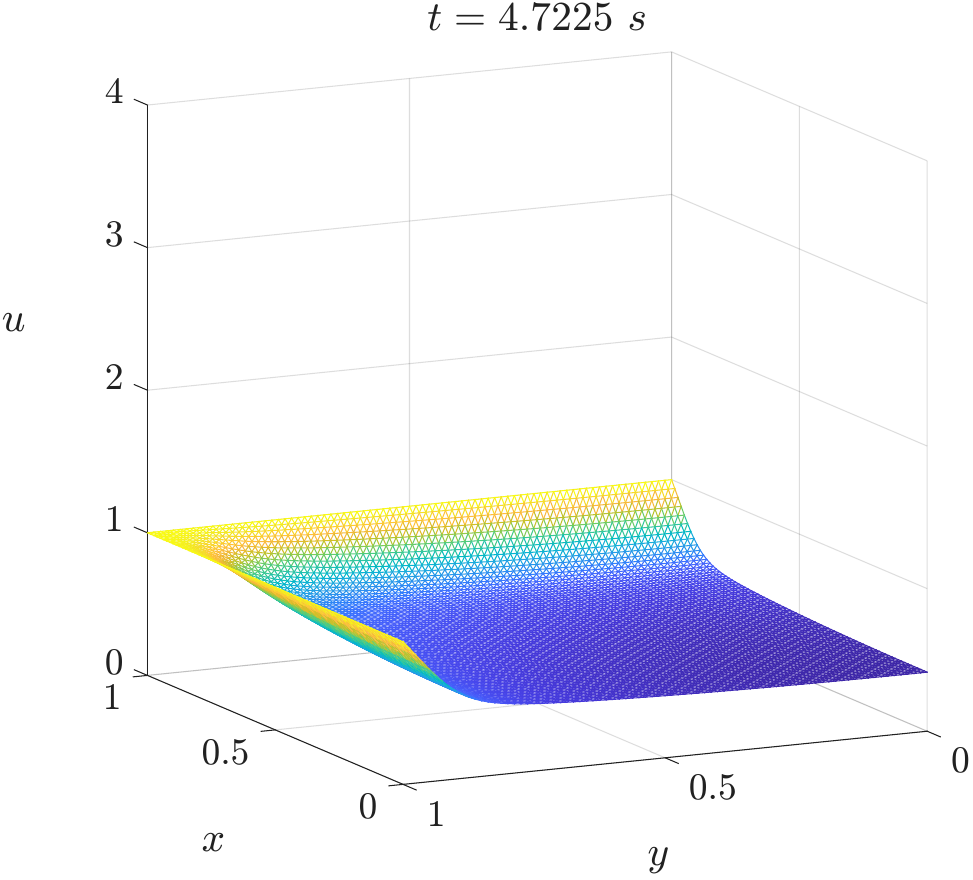}
		\includegraphics[width=0.24\linewidth]{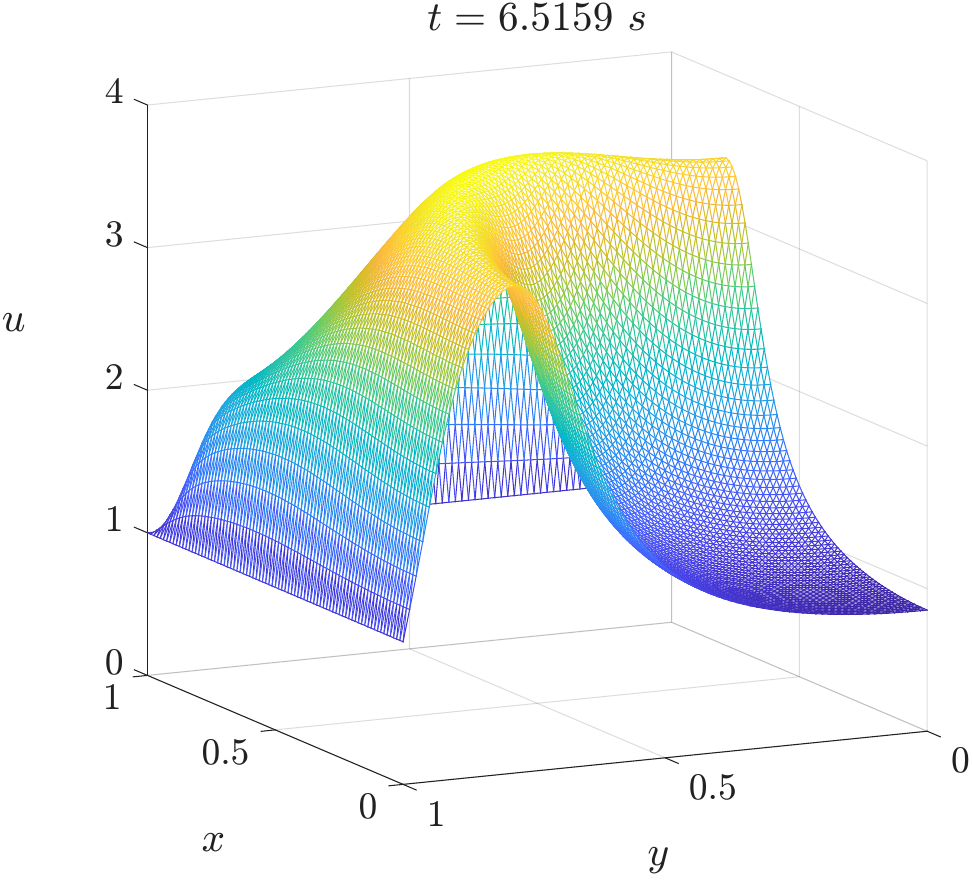}
	\end{minipage}
	
	\vspace{1cm} 
	
	\begin{minipage}{\textwidth}
		\centering
		\includegraphics[width=0.24\linewidth]{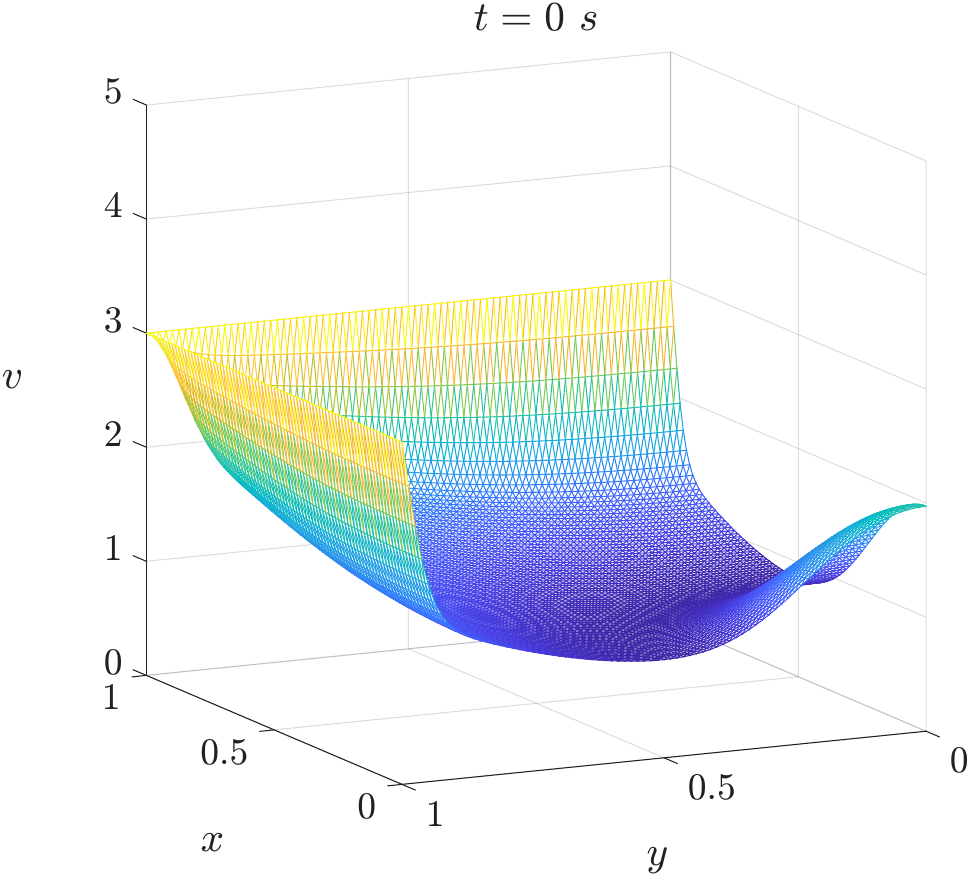}
		\includegraphics[width=0.24\linewidth]{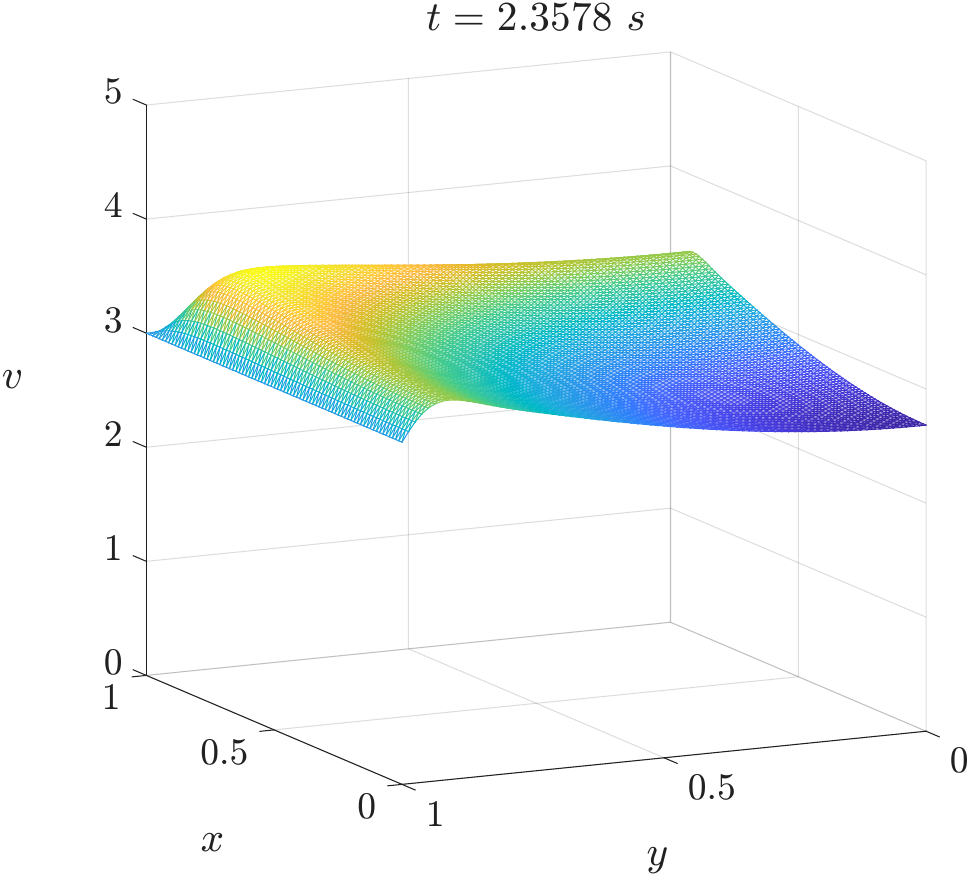} 
		\includegraphics[width=0.24\linewidth]{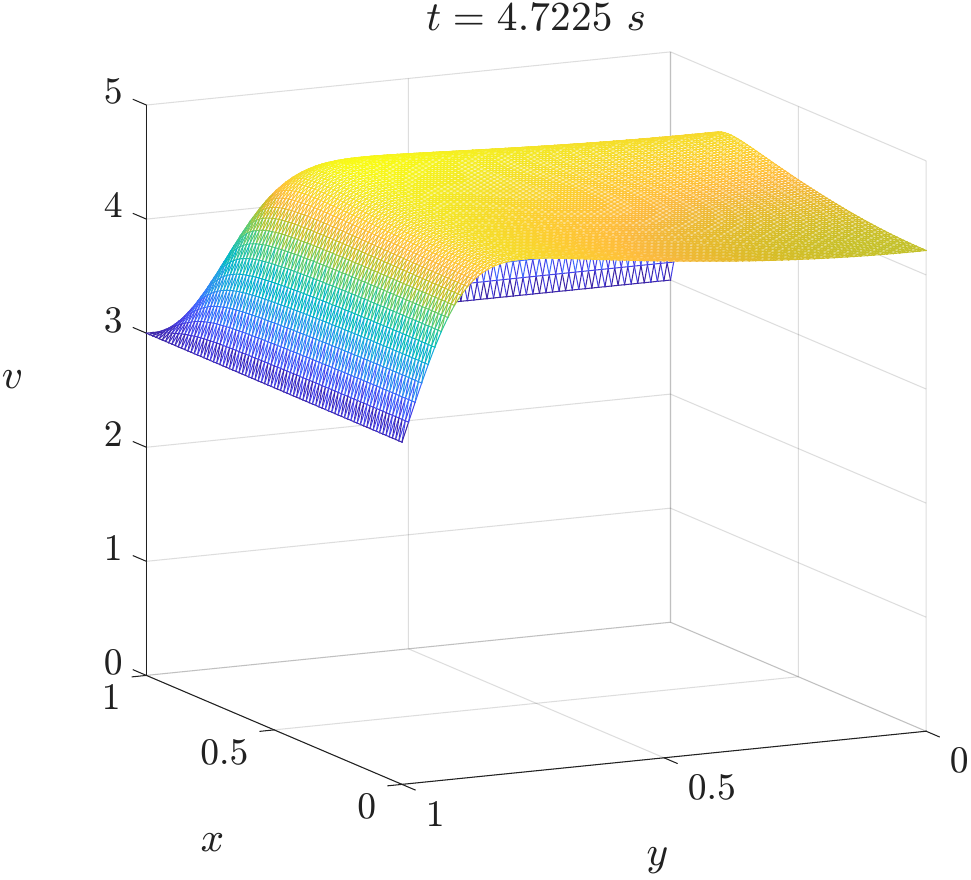}
		\includegraphics[width=0.24\linewidth]{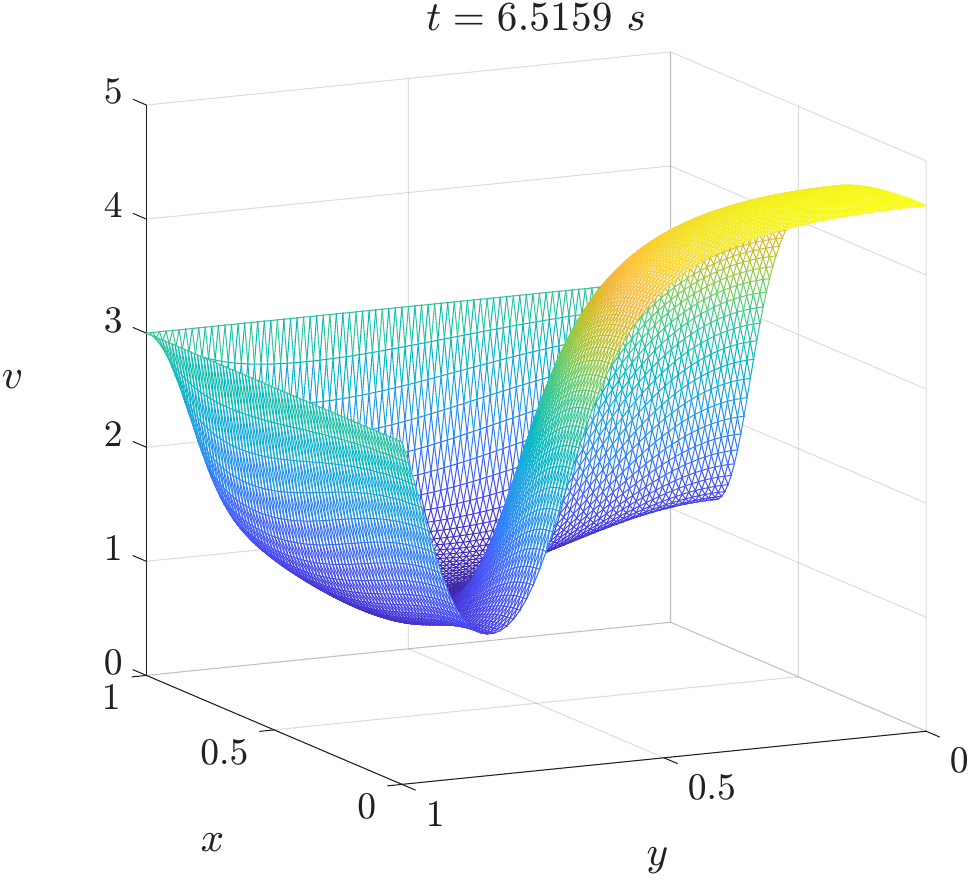}

	\end{minipage}
	\caption{Components $u$ (up) and $v$ (down) of the periodic solution at time $t=0, t=2.3578 \approx T/3, t=4.7225 \approx 2T/3$ and $t=6.5159$, where the maximum of $ \|\nabla u \|_0^2 + \|\nabla v\|_0^2 $  is attained. }
	\label{fig:femsolution}
\end{figure}

We choose the value $\nu=0.002$, which offered an adequate balance between spatial complexity and computability. In this case we consider a triangulation of~$\Omega$ based on a uniform mesh with 80  subdivision per side and with diagonals running southwest-northeast. We used a spatial discretization with quadratic finite elements, which resulted in a system with {\color{red}51200} degrees of freedom. For this discretization, we computed the corresponding limit cycle (see details in~\cite[Section~6]{temporal_nos}), the period being $T=7.090636$ (up to 7 significant digits). {\color{red}In Figure \ref{fig:femsolution}, we present the FEM solution at different significative time steps to see the complexity of the problem as time evolves. We consider $t=0$, $t\approx T/3$ and $t \approx 2T/3$. We also show the solution at $ t= 6.515906$ (up to 7 significant digits), where the maximum of $ \|\nabla u \|_0^2 + \|\nabla v\|_0^2 $  is attained, to see the large variations that occur at initial time and near the endpoints of the period.}

In all the experiments reported below the value of~$\tau$ in~(\ref{labU})  was set to 1. The mean is subtracted from the snapshots $\bu_h(t_j)$ to obtain the fluctuations. As a consequence, we transform the problem into one with homogeneous boundary conditions on~$\Gamma_1$ and $w_0$ in~(\ref{labU})  is taken as 0 (the mean value).

{\color{red}We consider the continuous-in-time problem:
find $u_r:(0,T]\rightarrow \bU^r$
such that
\begin{eqnarray}\label{eq:pod_ex}
	(u_{r,t},v_r)+\nu(\nabla u_r,\nabla v_r)+(g(u_r),v_r)=(f,v_r),\quad \forall\ v_r\in \bU^r,
\end{eqnarray}
with $u_r(0)=P^r u^0$. Then, we integrate in time (\ref{eq:pod_ex}) with the BDF-q $q=5$, $\Delta t =T/1024$, 
so that the temporal errors are negligible. In Tables \ref{tab:theorem4L2} and \ref{tab:theorem4H1}
we show the $L^2$ and $H^1$ norms of the errors: $\bu_r^n-P^r\bu_h(t_n)$ and $(I-P^r)\bu_h(t_n)$. We can observe that both errors decrease as $r$ increases and that the behavior of the error of the POD method from $r=34$ follows very closely that of the errors in the projection.}

\begin{table}[htbp]
	\centering
	\begin{tabular}{|c|c|c|}
		\hline
		$r$ & $\max\limits_{q\leq n \leq M}{\|\bu_r^n - P^r\bu_h(t_n) \|}_0$ &
		$\max\limits_{q\leq n \leq M}{\|(I - P^r) \bu_h(t_n)\|_0}$   \\
		\hline
		18 & $1.8020 \times 10^{-1}$ & $7.7239 \times 10^{-3}$ \\
		26 & $1.3197 \times 10^{-3}$ & $7.2341 \times 10^{-4}$ \\
		34 & $9.9256 \times 10^{-5}$ & $1.0012 \times 10^{-4}$ \\
		42 & $4.4413 \times 10^{-6}$ & $5.0967 \times 10^{-6}$ \\
		50 & $3.3172 \times 10^{-7}$ & $3.5240 \times 10^{-7}$ \\
		\hline
	\end{tabular}
	\caption{Maximum errors in $L^2$ norm for $q=5$.}
	\label{tab:theorem4L2}
\end{table}

\begin{table}[htbp]
	\centering
	\begin{tabular}{|c|c|c|}
		\hline
		$r$ & $\max\limits_{q\leq n \leq M}{\|\nabla(\bu_r^n - P^r\bu_h(t_n))\|}_0$ &
		$\max\limits_{q\leq n \leq M}{\|\nabla(I - P^r) \bu_h(t_n)\|_0}$   \\
		\hline
			18 & $9.0212 \times 10^{-1}$ &  $8.7820 \times 10^{-2}$ \\
			26 & $1.1484 \times 10^{-2}$ &  $8.7337 \times 10^{-3}$ \\
			34 & $6.6843 \times 10^{-4}$ &  $7.2388 \times 10^{-4}$ \\
 			42 & $6.1572 \times 10^{-5}$ &  $6.0239 \times 10^{-5}$ \\
			50 & $1.0849 \times 10^{-5}$ &  $4.4051 \times 10^{-6}$ \\
		\hline
	\end{tabular}
	\caption{Maximum errors in $H_0^1$ norm for $q=5$.}
	\label{tab:theorem4H1}
\end{table}

{\color{red}To check the advantage of high order time integrators in the numerical integration of POD methods, in Table \ref{tab:timeqsL2} we show the maximum error of the POD approximation for different values of $r$ and $q$. In this table we keep $\Delta t =T/1024$. It can be observed that, except for first value of $r$, second order is necessary for optimal errors for $r=26$, third order for $r=32$, forth order for $r=42$ and fifth order for $r=50$. This means that as $r$ increases the error is determined by the time integrator and higher order methods allow to achieve optimal errors without the need to increase the number of time steps.}

\begin{table}[htbp]
	\centering
	\begin{tabular}{|c|ccccc|}
		\hline
		$r$ 
		& \multicolumn{5}{|c|}{$\max\limits_{q\leq n \leq M} \|\bu_r^n - \bu_h(t_n)\| $ for different $q$} \\
		\cline{2-6}
		 & $q=1$ & $q=2$ & $q=3$ & $q=4$ & $q=5$ \\
		\hline
		18 & $1.2915 \times 10^{-1}$ & $1.8004 \times 10^{-1}$  & $1.7983 \times 10^{-1}$ &  $1.7983 \times 10^{-1}$ &  $1.7983 \times 10^{-1}$\\
		26 & $5.3288 \times 10^{-2}$ & $1.5079\times 10^{-3}$ & $1.3063 \times 10^{-3}$ &  $1.3020 \times 10^{-3}$ &  $1.3020\times 10^{-3}$\\
		32 & $4.2505 \times 10^{-2}$ & $7.9057\times 10^{-4}$ &   $6.5009\times 10^{-5}$ &  $5.9256\times 10^{-5}$ &  $5.9739\times 10^{-5}$\\
		42 &  $3.8688 \times 10^{-2}$ &  $7.7619 \times 10^{-4}$ &  $2.6741\times 10^{-5}$ &  $3.1813\times 10^{-6}$ &  $2.2734 \times 10^{-6}$ \\
		50 &  $6.3663 \times 10^{-2}$ &  $1.0528\times 10^{-3}$ &   $2.9448 \times 10^{-5}$ &  $2.4092 \times 10^{-6}$  & $3.0961 \times 10^{-7}$\\
		\hline
	\end{tabular}
	\caption{Maximum errors in $L^2$ norm for $1\le q\le 5$ and different values of $r$.}
	\label{tab:timeqsL2}
\end{table}

{\color{red}To check numerically the rate of convergence in time, we integrate (\ref{eq:pod_ex})} with a numerical differentiation formulae (NDF) \cite{NDF1}, 
in the variable-step, variable-order implementation of {\sc Matlab}'s command {\tt ode15s} \cite{odesuite}, with sufficiently small tolerances for the local errors ($10^{-12}$). We are interested in checking the theoretical rate of convergence in time of the POD method based on the BDF-q formula. We fix $r=18$ in (\ref{eq:bU_r}) and take different values of $\Delta t =T/M$. 
Then, in our numerical experiments we measure
$u_{r}^{s}-u_{r}^{bdf_q}$, for $u_r^{bdf_q}$ the solution of (\ref{eq:pod_ex}) with BDF-q. We observe that with this procedure the only component of the
error in (\ref{er_defi_le3}) is the second term on the right-hand side:
\begin{equation}\label{eq:rhs2term}
	C e^{2\alpha T}(\Delta t)^{2q}\int_0^T\|\partial_t^{q+1}u_h(s)\|_0^2 \ ds.
\end{equation}

\begin{figure}[!]
	\centering
	\begin{minipage}[b]{0.45\linewidth}
		\centering
		\includegraphics[width=\linewidth]{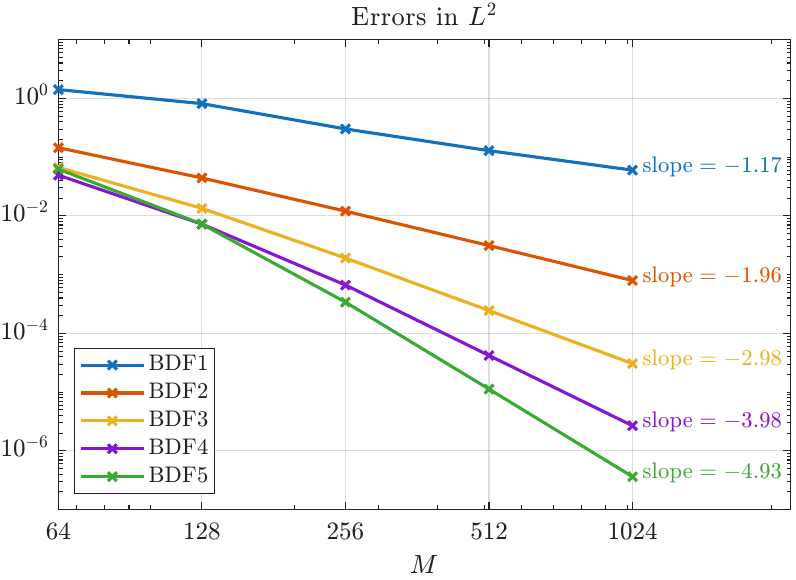}
	\end{minipage}
	\hfill
	\begin{minipage}[b]{0.45\linewidth}
		\centering
		\includegraphics[width=\linewidth]{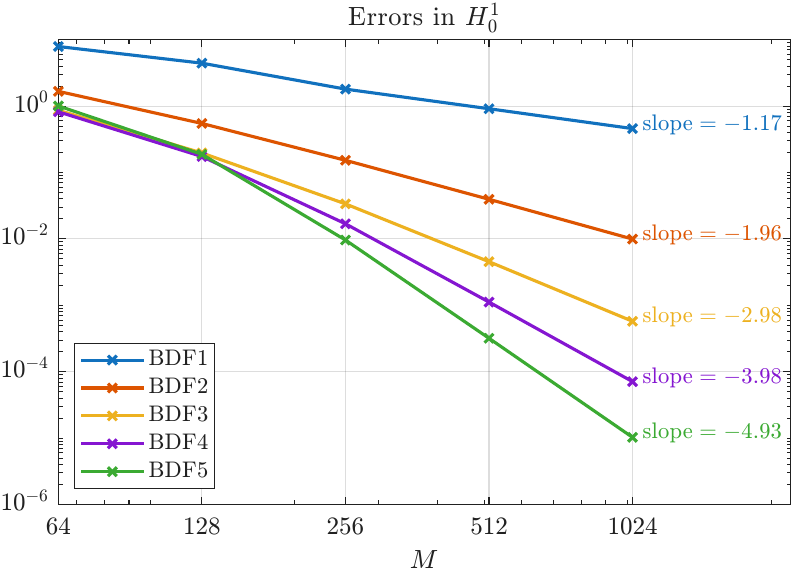}
	\end{minipage}
	\caption{Maximum error along one period in the $L^2$ norm, $\max\limits_{q\leq n \leq M}\|u_r^s(t_n) - u_r^{bdf_q}(t_n)\|_0,$ (left) and
		$H_0^1$ norm, $\max\limits_{q\leq n \leq M}\|\nabla(u_r^s(t_n) - u_r^{bdf_q}(t_n))\|_0,$ (right) for $M=64, 128, 256, 512$ and $1024$.}
	
	\label{fig:maxerr}
\end{figure}

To compute the starting values in $u_{r}^{bdf_q}$ we use the following procedure. For $q=1$ the only starting value $u_{r}^{bdf_1}(t_0)=P^r u^0$. For $q=2$,
$u_{r}^{bdf_2}(t_0)=P^r u^0$ as before. To compute $u_{r}^{bdf_2}(t_1)$ we solve (\ref{eq:pod_ex}) from $t_0$ to $t_1$ with the BDF-1 method.
For $q\ge 3$ the BDF-q  computes the starting values at times $t_1$ to $t_{q-1}$ with the $q-1$ formula and fixed step $\Delta t'=(\Delta t )^{q/(q-1)}$
to guarantee that the  starting values have order $q$ in time. In tables \ref{table1} and \ref{table2} we show the maximum  error of the starting values of the formulas of order $2$ to $5$ in the $L^2$ and $H_0^1$ norms respectively, for $M=64, 128, 256, 512$ and $1024$, i.e. $\max\limits_{1\leq n \leq q-1}\|u_r^s(t_n) - u_r^{bdf_q}(t_n) \|_0$ and
$\max\limits_{1\leq n \leq q-1}\|\nabla(u_r^s(t_n) - u_r^{bdf_q}(t_n)) \|_0$, respectively. The same values of $M$ are used in Figure \ref{fig:maxerr}, where 
we plot the maximum values of the error along the full period in both norms. We can observe that all the errors in Tables \ref{table1} and \ref{table2} are well below the
maximum errors of Figure \ref{fig:maxerr} which means that the described procedure to choose the starting values has essentially no influence in the error of the method.

\begin{table}[h!]
	\centering
	\begin{tabular}{|c|c|c|c|c|}
		\hline
		& $q=2$ & $q=3$ & $q=4$ & $q=5$ \\
		\hline
		$M = 64$ & \(5.647 \times 10^{-2}\) & \(5.138 \times 10^{-3}\) & \(7.707 \times 10^{-4}\) & \(1.245 \times 10^{-3}\) \\
		$M = 128$ & \(2.463\times 10^{-2}\) & \(1.522 \times 10^{-3}\) & \(1.921 \times 10^{-4}\) & \(5.928 \times 10^{-6}\) \\
		$M = 256$ & \(7.910 \times 10^{-3}\) & \(2.461 \times 10^{-4}\) & \(1.292 \times 10^{-5}\) & \(4.172 \times 10^{-6}\) \\
		$M = 512$ & \(2.212 \times 10^{-3}\) & \(4.129 \times 10^{-5}\) & \(5.688 \times 10^{-7}\) & \(1.865 \times 10^{-6}\) \\
		$M = 1024$ & \(5.830 \times 10^{-4}\) & \(5.191 \times 10^{-6}\) & \(1.999 \times 10^{-9}\) & \(2.679 \times 10^{-10}\) \\
		\hline
	\end{tabular}
	\caption{\label{table1} 
		Errors at the starting values in the $L^2$ norm: $\max\limits_{1\leq n \leq q-1}\|u_r^s(t_n) - u_r^{bdf_q}(t_n) \|_0.$
	}
\end{table}

\begin{table}[h!]
	\centering
	\begin{tabular}{|c|c|c|c|c|}
		\hline
		& $q=2$ & $q=3$ & $q=4$ & $q=5$ \\
		\hline
		$M = 64$  & \(2.911 \times 10^{-1}\) & \(2.655 \times 10^{-2}\) & \(7.860 \times 10^{-3}\) & \(1.395 \times 10^{-2}\) \\
		$M = 128$ & \(1.214 \times 10^{-1}\) & \(7.524 \times 10^{-3}\) & \(1.893 \times 10^{-3}\) & \(4.806 \times 10^{-4}\) \\
		$M = 256$ & \(3.954\times 10^{-2}\) & \(1.248 \times 10^{-3}\) & \(1.217 \times 10^{-4}\) & \(3.920 \times 10^{-5}\) \\
		$M = 512$ & \(1.144 \times 10^{-2}\) & \(2.160 \times 10^{-4}\) & \(5.660 \times 10^{-6}\) & \(2.724 \times 10^{-6}\) \\
		$M = 1024$ & \(3.090\times 10^{-3}\) & \(2.770 \times 10^{-5}\) & \(2.389 \times 10^{-8}\) & \(4.049 \times 10^{-9}\) \\
		\hline
	\end{tabular}
	\caption{\label{table2} 
		Errors at the starting values in the $H_0^1$ norm: $\max\limits_{1\leq n \leq q-1}\|\nabla(u_r^s(t_n) - u_r^{bdf_q}(t_n)) \|_0.$
	}
\end{table}

Since the BDF-q formulas are implicit we use Newton's method with standard local extrapolation for the initial guess  to solve the corresponding equations. 
We take tolerance  $(\Delta t )^q /100$ to guaranty that the error coming from the Newton's iteration does not affect the rate of convergence of the corresponding BDF formula. We have observed that the number of steps in the Newton's iteration ranges between 2 and 4. The average number of steps in Newton's method across one period considering all methods is $2.97$.
As stated before, in Figure \ref{fig:maxerr} we represent the maximum error in time in one period both in $L^2$ and $H^1$ norms. 

The expected rate of convergence $q$ can be clearly observed for all the formulas, according  to (\ref{eq:rhs2term}). We have measured the slopes
in the figure in the asymptotic regime. While for formulas 1 and 2 the asymptotic slope can be observed from the first value of $M$, $M=64$, for the last two formulas 4 and 5 this regime needs a larger value of $M$. This fact is in agreement with our theoretical results since condition (\ref{cond_deltat1})
is required to achieve the rate of convergence. This means that smaller $\Delta t$ is required as $q$ grows since $\eta_q$ increases in the
right condition of (\ref{cond_deltat1}) while $\tilde\lambda_m$ decreases.
\bigskip
{\color{red}\section{Conclusions}

In this paper we have obtained a priori bounds for fully discrete POD methods applied to reaction-diffusion equations
using BDF-$q$ schemes, $q\le 5$. Error bounds for POD methods with high-order time discretizations cannot be found in the literature. To get those bounds, we follow reference \cite{bdf_stokes} and apply $G$-stability techniques to handle the temporal error. The analysis of a reduced-order method also requires an error estimate of the projection of the approximation to the time derivative given by the BDF-q method onto the reduced-order space. To bound this term we use first-order finite differences in the set of snapshots and prove that the approximation to the time derivative in
BDF-$q$ method can always be written as a linear combination of those first-order differences. The numerical experiments confirm the rate of convergence predicted by the theory.
}

{\bf Data Availability} The authors declare that the data supporting the findings of this study can be obtained by reproducing the computations as described in the section on numerical experiment. The data and the MATLAB source code to reproduce the numerical experiments will be located in a public repository when the paper is accepted for publication.
\bigskip

{\bf Declarations}
Competing interests. The authors declare that they have no conflict of interest.



\begin{thebibliography}{10}
\bibitem{il_sur} F. Ballarin and T. Iliescu.
\newblock{A Priori Error Bounds for POD-ROMs for Fluids: A Brief Survey.}
{\rm In: Giacomini, M., Perotto, S., Rozza, G. (eds) Emerging Technologies in Computational Sciences for Industry, Sustainability and Innovation. M2P 2023. Lecture Notes in Computational Science and Engineering, vol 146. Springer, Cham.}

\bibitem{brenner-scot}
S.~C. Brenner and L.~R. Scott.
\newblock {\em The mathematical theory of finite element methods}, volume~15 of
  {\em Texts in Applied Mathematics}.
\newblock Springer, New York, third edition, 2008.

\bibitem{bdf_stokes} A. Contri, B. Kov\'acs and A. Massing.
\newblock Error analysis of BDF 1-6 time-stepping methods for the transient stokes problems: velocity and pressure estimates.
\newblock {\em SIAM J. Numer. Anal.},  (to appear). 


\bibitem{dalquis1} G. Dahlquist.
\newblock{Error analysis for a class of methods for stiff non-linear initial value problems.}
\newblock Lecture Notes in Mathematics, Springer Berlin Heidelberg, 60-72, 1976.

\bibitem{dalquis2} G. Dahlquist.
\newblock{G-stability is equivalent to A-stability.}
\newblock {\em BIT}, 18, 384-401, 1978.


\bibitem{letter} B.~Garc\'{\i}a-Archilla, V. John and J. Novo.
\newblock Second order error bounds for POD-ROM methods based on first order divided differences
\newblock {\em Appl. Math. Lett.}, 146, Paper No. 108836, 7 pp. 2023.

\bibitem{temporal_nos} B.~Garc\'{\i}a-Archilla, V. John and J. Novo.
\newblock POD-ROM methods: from a finite
set of snapshots to continuous-in-time approximations
\newblock {\em SIAM J. Numer. Anal.} 63, 800-826, 2025.

\bibitem{nos_bdf} B.~Garc\'{\i}a-Archilla and J. Novo.
\newblock Robust error bounds for the Navier-Stokes equations using implicit-explicit second order 
BDF method with variable steps
\newblock {\em IMA J. Numer. Anal.} 43, 2892-2933, 2023.                                                   

\bibitem{Bosco_pointwise} B.~Garc\'{\i}a-Archilla and J. Novo.
\newblock Pointwise error bounds in POD methods without difference quotients
\newblock {\em J. Sci. Comput,} 103, No. 1, paper 24, 26 pp, 2025.

\bibitem{bosco-titi-fem}
B.~Garc\'{\i}a-Archilla and E.~S. Titi.
\newblock Postprocessing the {G}alerkin method: the finite-element case.
\newblock {\em SIAM J. Numer. Anal.}, 37, 470--499, 2000.


\bibitem{Hairer}
E.~Hairer, S.~P. N{\o}rsett, and G.~Wanner.
\newblock {\em Solving ordinary differential equations. {I}}, volume~8 of {\em
  Springer Series in Computational Mathematics}.
\newblock Springer-Verlag, Berlin, second edition, 1993.
\newblock Nonstiff problems.


\bibitem{NDF1}
R.~W. Klopfenstein.
\newblock Numerical differentiation formulas for stiff systems of ordinary
  differential equations.
\newblock {\em RCA Rev.}, 32:447--462, 1971.


\bibitem{rubino_etal} B. Koc, S. Rubino, M. Schneier, J. Singler and T. Iliescu.
\newblock On optimal pointwise in time error bounds and difference quotients for the proper orthogonal decomposition.
\newblock {\em SIAM J. Numer. Anal.},  59, 2163-2196, 2021.

\bibitem{ku-Vol} K. Kunisch and S. Volkwein.
\newblock Galerkin proper orthogonal decomposition methods for parabolic problems.
\newblock {\em Numer. Math.}, 90, 117-148, 2001.


\bibitem{nevan_ode} O. Nevanlinna and F. Odeh.
\newblock Multiplier techniques for linear multistep methods.
\newblock{\em Numer. Funct. Anal. and Opt.}, 3, 377--423, 1981.



\bibitem{odesuite}
L.~F. Shampine and M.~W. Reichelt.
\newblock The {MATLAB} {ODE} suite.
\newblock {\em SIAM J. Sci. Comput.}, 18(1):1--22, 1997. 
\newblock Dedicated to C. William Gear on the occasion of his 60th birthday.


\bibitem{Thomee} V. Thom\'ee.
\newblock{\em Galerkin finite elements for parabolic problems}.
\newblock Springer Verlag, Heidelberg, Germany, 2006.

\bibitem{Wang} Z. Wang.
\newblock Nonlinear Model Reduction Based on the Finite Element Method with Interpolated Coefficients: Semilinear Parabolic 
Equations.
\newblock{\em Numer. Methods Partial Differential Equations}, 31, 1713--1741, 2015.

\end{thebibliography}
\end{document}